\documentclass[preprint,12pt]{elsarticle}

\usepackage{lipsum}
\usepackage{amsfonts}
\usepackage{mathtools}
\usepackage{amsopn}
\usepackage{graphicx}
\usepackage{tikz}
\usetikzlibrary{trees,calc}
\usepackage{epstopdf}
\usepackage{algorithmic}
\usepackage{multirow}
\usepackage{booktabs}
\usepackage{amsmath}
\usepackage{color}
\usepackage{upgreek}
\usepackage{amssymb}
\usepackage{mathrsfs}
\usepackage{bm}
\usepackage{todonotes}
\ifpdf
  \DeclareGraphicsExtensions{.eps,.pdf,.png,.jpg}
\else
  \DeclareGraphicsExtensions{.eps}
\fi

\usepackage{xcolor} 
\usepackage{graphicx}
\usepackage{subcaption}
\graphicspath{{Pictures/}}

\newcommand{\n}{\mathbf{n}}
\newcommand{\x}{\mathbf{x}}

\allowdisplaybreaks 

\usepackage{amsthm} 
\newtheorem{remark}{Remark}

\usepackage[hidelinks]{hyperref}
\usepackage{cleveref}

\begin{document}

\begin{frontmatter}

\title{Augmented Lagrangian preconditioning for a simplified Ericksen--Leslie model of nematic liquid crystals}

\author{Yanying Li} 
\author{Xu Qian}
\author{Jingmin Xia\corref{cor1}}
\affiliation{organization={College of Science, National University of Defense Technology},
            addressline={liyanying@nudt.edu.cn}, 
            city={Changsha},
            country={China}}

\affiliation{organization={College of Science, National University of Defense Technology},
            addressline={qianxu@nudt.edu.cn}, 
            city={Changsha},
            country={China}}

\affiliation{organization={College of Meteorology and Oceanography, National University of Defense Technology},
            addressline={jingmin.xia@nudt.edu.cn}, 
            city={Changsha},
            country={China}}
\cortext[cor1]{Corresponding author.}

\begin{abstract}
The numerical solution of the simplified Ericksen--Leslie model for nematic liquid crystals is challenging because the flow and director equations are strongly coupled and because incompressibility and the unit-length condition must be enforced simultaneously. A Lagrange multiplier formulation avoids a small Ginzburg--Landau parameter, but the Newton systems have a double saddle-point structure. We develop an augmented Lagrangian block preconditioner in which both constraints are augmented while their discrete enforcement remains multiplier based. After finite element discretization and backward Euler time integration, the Newton increments are grouped into velocity--director and pressure--multiplier variables. A block-diagonal approximation of the coupled velocity--director block then leads to separate, physically scaled approximations of the pressure and director-multiplier Schur complements. Manufactured-solution tests show the expected spatial accuracy and first-order temporal convergence for the primary variables; the multiplier error reaches a spatial-error floor on the fixed mesh used in the temporal study. In the reported parameter ranges, the outer FGMRES iteration counts are nearly mesh independent, remain stable under time-step and viscosity variation, and improve as the augmentation parameters increase. A smooth benchmark also exhibits monotone decay of the computed total energy.
\end{abstract}


\end{frontmatter}

\section{Introduction}

Nematic liquid crystals are a class of soft matter that flows like a viscous fluid while exhibiting long-range orientational order. 
The local molecular alignment is described by a unit vector field, the director \(\mathbf{n}(\mathbf{x},t)\), 
which satisfies \(|\mathbf{n}|=1\) pointwise. This combination of fluidity and anisotropy underpins numerous technological applications 
and also leads to rich phenomena such as the formation, motion and annihilation of topological defects\cite{deGennesProst1993,Stewart2004}.

The hydrodynamic theory of nematic liquid crystals was formulated by Ericksen \cite{ericksen1961, ericksen1962ARMA} 
and Leslie \cite{leslie1968ARMA, leslie1979}, who proposed a model incorporating the evolution of the director field coupled 
with the incompressible Navier--Stokes equations through an additional elastic stress tensor arising from molecular alignment. 
To simplify the analysis while retaining the essential mathematical structure, Lin \cite{lin1989,LinLiu1995} considered the regime 
in which the Leslie viscous stresses are neglected except for an isotropic viscosity, and combined with the one-constant approximation 
of the Oseen--Frank elastic energy \cite{Oseen1933,Frank1958}, introduced a simplified Ericksen--Leslie model, which reads in dimensionless form
\begin{equation}\label{eq:EL-original}
\begin{cases}
\mathbf{n}_t + \mathbf{u}\cdot\nabla\mathbf{n} - \mu\Delta\mathbf{n} - \mu|\nabla\mathbf{n}|^2\mathbf{n} = 0,\\[4pt]
\mathbf{u}_t + \mathbf{u}\cdot\nabla\mathbf{u} + \nabla p - \nu\Delta\mathbf{u} + K\nabla\cdot\big((\nabla\mathbf{n})^\top\nabla\mathbf{n}\big) = 0,\\[4pt]
\nabla\cdot\mathbf{u} = 0,\\[4pt]
|\mathbf{n}| = 1,
\end{cases}
\end{equation}
with the boundary condition
\begin{equation}
    \mathbf{u}|_{\partial\Omega} = 0,\quad \partial_{\boldsymbol{\nu}}\n|_{\partial\Omega} = 0,
\end{equation}
and initial conditions
\begin{equation}
    \mathbf{u}(\x,0) = \mathbf{u}_0,\qquad \n(\x,0) = \n_0, \qquad \x\in\Omega\subset \mathbb{R}^{d}.
\end{equation}
Here \(\mathbf{n}\) is the director, \(\mathbf{u}\) is the velocity, \(p\) is the pressure, and the positive constants \(\nu\), \(\mu\), and \(K\) 
denote the fluid viscosity, the director relaxation (mobility) coefficient, and the Frank elastic constant, respectively. The system is supplemented with no-slip boundary conditions 
for \(\mathbf{u}\) and homogeneous Neumann conditions for \(\mathbf{n}\), together with initial data satisfying \(|\mathbf{n}_0|=1\) and \(\nabla\cdot\mathbf{u}_0=0\).

System \eqref{eq:EL-original} possesses an intrinsic energy dissipation law~\cite{LinLiu1995}. Defining the total free energy
\[
E(\mathbf{u},\mathbf{n}) = \frac12\|\mathbf{u}\|^2 + \frac{K}{2}\|\nabla\mathbf{n}\|^2,
\]
where \(\|\cdot\|\) is the \(L^2(\Omega)\) norm, smooth solutions satisfy
\begin{equation}\label{eq:energy-law-orig}
\frac{\mathrm{d}}{\mathrm{d}t}E(\mathbf{u},\mathbf{n}) = -\nu\|\nabla\mathbf{u}\|^2 - K\mu\|\Delta\mathbf{n} + |\nabla\mathbf{n}|^2\mathbf{n}\|^2 \le 0.
\end{equation}
This dissipative structure underlies the global existence of weak solutions and local well-posedness of classical solutions established by Lin and Liu \cite{LinLiu1995}.

The main numerical challenge of \eqref{eq:EL-original} stems from the non-convex unit-length constraint \(|\mathbf{n}|=1\). Conventional strategies for handling such constraints---including projection methods \cite{Alouges1997}, Ginzburg--Landau penalty approaches \cite{LinLiu1995}, and Lagrange multiplier techniques---each have known limitations. Projection may compromise temporal accuracy or energy stability, penalty methods lead to ill-conditioned algebraic systems as the penalty parameter \(\epsilon\to0\), and the Lagrange multiplier approach introduces a saddle-point structure that requires inf--sup compatible discretizations.

A particularly effective Lagrange multiplier formulation for this problem was proposed by Badia \cite{Badia2011}, who introduced a scalar field \(q(\mathbf{x},t)\) to enforce the unit-length constraint explicitly. Using the identity \(\mathbf{n}\cdot\Delta\mathbf{n} = -|\nabla\mathbf{n}|^2\), which holds when \(|\mathbf{n}|=1\), the cubic term \(-|\nabla\mathbf{n}|^2\mathbf{n}\) in \eqref{eq:EL-original} is replaced by \(q\mathbf{n}\), yielding the equivalent saddle-point system
\begin{equation}\label{eq:EL-saddle}
\begin{cases}
\mathbf{n}_t + \mathbf{u}\cdot\nabla\mathbf{n} + \mu(-\Delta\mathbf{n} +  q\,\mathbf{n}) = 0,\\[4pt]
\mathbf{u}_t + \mathbf{u}\cdot\nabla\mathbf{u} + \nabla p - \nu\Delta\mathbf{u} + K\nabla\cdot\big((\nabla\mathbf{n})^\top\nabla\mathbf{n}\big) = 0,\\[4pt]
\nabla\cdot\mathbf{u} = 0,\\[4pt]
|\mathbf{n}| = 1.
\end{cases}
\end{equation}

\begin{remark}
When \(|\mathbf{n}|=1\), the identities \(\mathbf{n}\cdot\mathbf{n}_t=0\), \(\mathbf{n}\cdot(\mathbf{u}\cdot\nabla\mathbf{n})=0\), and
\(
\tfrac12\Delta|\mathbf{n}|^2 = \mathbf{n}\cdot\Delta\mathbf{n} + |\nabla\mathbf{n}|^2
\)
imply \(-\mathbf{n}\cdot\Delta\mathbf{n}=|\nabla\mathbf{n}|^2\). Taking the Euclidean inner product of the director equation in \eqref{eq:EL-saddle} with \(\mathbf{n}\) therefore gives \(q=-|\nabla\mathbf{n}|^2\), which recovers the cubic term in \eqref{eq:EL-original}. Thus the two systems are equivalent for sufficiently regular solutions satisfying the unit-length constraint.
\end{remark}

\begin{remark}
The elastic stress tensor in the momentum equation can be rewritten using the identity \cite{Badia2011,CaoYi2025}
\(
\nabla\cdot\big((\nabla\mathbf{n})^\top\nabla\mathbf{n}\big) = (\nabla\mathbf{n})^\top\Delta\mathbf{n} + \frac12\nabla(|\nabla\mathbf{n}|^2).
\) 
The second term is the gradient of a scalar field and can be absorbed into the pressure by defining the modified pressure \(\widetilde{p} = p + \frac{K}{2}|\nabla\mathbf{n}|^2\). The momentum equation then becomes
\(
\mathbf{u}_t + \mathbf{u}\cdot\nabla\mathbf{u} - \nu\Delta\mathbf{u} + \nabla\widetilde{p} + K(\nabla\mathbf{n})^\top\Delta\mathbf{n} = 0.
\)
\end{remark}

At the continuous level one has \(q = -|\nabla\mathbf{n}|^2\) under the sign convention in \eqref{eq:EL-saddle}, 
and the energy law retains its form. This reformulation reduces the cubic nonlinearity to a bilinear coupling \(q\mathbf{n}\), 
simplifies Newton linearization, and, importantly, makes the stability of the multiplier depend on a discrete inf-sup condition rather 
than on a penalty parameter; consequently, conditioning and solvability become questions of discrete compatibility and block solver design 
even when the constraint is enforced exactly \cite{Badia2011}. The broader saddle-point framework of Badia et al.~\cite{Badia2011} also treats Ginzburg--Landau penalized and exactly constrained variants in a unified formulation.

Building on this formulation, Cao and Yi \cite{CaoYi2025} developed fully decoupled, length-preserving predictor-corrector schemes with second-order temporal accuracy, and later combined a pressure-correction strategy with the scalar auxiliary variable  approach to achieve linear, unconditionally stable schemes for the original simplified model \eqref{eq:EL-original} \cite{cao2026}. Despite these advances, solving the fully discrete systems arising from \eqref{eq:EL-saddle} in a monolithic fashion remains computationally demanding: at each time step the linearized problem is a \(4\times4\) block matrix with a double saddle-point structure coupling \((\mathbf{u},p)\) on one side and \((\mathbf{n},q)\) on the other. The design of robust and scalable iterative solvers for such systems is an active area of research.

The augmented Lagrangian (AL) method combines the advantages of penalty and multiplier approaches: the Lagrange multiplier preserves accuracy with moderate penalty parameters, while the penalty term improves the spectral properties of the Schur complement, facilitating the iterative solution of the linearized systems. Recently, Xia, Farrell and Wechsung \cite{xia2021} demonstrated that the AL strategy can effectively overcome the ill-conditioning associated with the unit-length constraint in the static Oseen--Frank model. In the present work, we extend this methodology to the fully time-dependent Ericksen--Leslie equations, applying the AL technique simultaneously to the incompressibility condition and the unit-length constraint.

The main contributions are threefold. First, we formulate a fully discrete augmented Lagrangian system for the time-dependent simplified Ericksen--Leslie equations, with both constraints retained in multiplier form. Second, we derive a block preconditioner for the Newton systems by separating the coupled velocity--director block from the pressure--multiplier block and constructing physically scaled approximations of the two constraint Schur complements. Third, we assess accuracy and solver behavior using manufactured solutions and benchmark computations covering mesh refinement, time-step and viscosity variation, a curved domain, and a two-defect configuration.

The remainder of the paper is organized as follows. In Section~\ref{sec:Augmented} we derive the augmented Lagrangian weak form. Section~\ref{sec:Finite} describes spatial and temporal discretizations and the Newton linearization. Preconditioning techniques are presented in Section~\ref{sec:Pre}. Numerical experiments are reported in Section~\ref{sec:Numerical}, and conclusions are drawn in Section~\ref{sec:conclusions}.

\section{Augmented Lagrangian methods}
\label{sec:Augmented}

The simplified Ericksen--Leslie model contains two equality constraints: the incompressibility constraint
\(\nabla\cdot\mathbf{u}=0\) and the pointwise unit-length constraint \(|\mathbf{n}|=1\). Direct discretization of these constraints leads to 
a coupled saddle-point algebraic system. The augmented Lagrangian method adds quadratic terms associated with the constraints while retaining the 
Lagrange multipliers; it therefore improves the conditioning of the constraint Schur complements without replacing the exact weak enforcement by a pure penalty method. 

For the present time-dependent dissipative system, the augmented Lagrangian functional introduced below should be understood as a device for identifying the multiplier and augmentation contributions to the weak equations, rather than as a complete variational principle for the full convective dynamics. We first describe the constraint part of the augmented functional and then state the corresponding augmented weak form used in the discretization and preconditioner construction.

\subsection{Establishment of the augmented Lagrangian functional}
Let \(\mathcal{J}(\mathbf{u},\mathbf{n})\) denote the energy part associated with the simplified Ericksen--Leslie model \eqref{eq:EL-original} on 
a bounded domain \(\Omega\subset\mathbb{R}^d\). We introduce the pressure \(p\) as the multiplier for \(\nabla\cdot\mathbf{u}=0\), and a scalar 
field \(q\) as the multiplier for \(|\mathbf{n}|^2-1=0\).  The constraint part of the Lagrangian is written as
\begin{equation}
\mathcal{L}(\mathbf{u},\mathbf{n},p,q)
=
\mathcal{J}(\mathbf{u},\mathbf{n})
-\left(p,\nabla\cdot\mathbf{u}\right)
+\frac{\mu}{2}\left(q,|\mathbf{n}|^2-1\right),
\end{equation}
where \((\cdot,\cdot)\) denotes the standard \(L^2(\Omega)\) inner product.  
The factor \(\mu/2\) in the director constraint is immaterial for the constraint equation, 
but it is convenient because variation with respect to \(\mathbf{n}\) gives the term \(\mu q\mathbf{n}\) in \eqref{eq:EL-saddle}.

We then introduce two positive augmentation parameters \(\gamma_u>0\) and \(\gamma_n>0\), associated respectively with 
incompressibility and the director-length constraint, and define
\begin{equation}
\begin{aligned}
\mathcal{L}_A(\mathbf{u},\mathbf{n},p,q;\gamma_u,\gamma_n)
={}&\mathcal{J}(\mathbf{u},\mathbf{n})
-\left(p,\nabla\cdot\mathbf{u}\right)
+\frac{\gamma_u}{2}\|\nabla\cdot\mathbf{u}\|_{L^2(\Omega)}^2 \\
&+\frac{\mu}{2}\left(q,|\mathbf{n}|^2-1\right)
+\frac{\gamma_n}{2}\||\mathbf{n}|^2-1\|_{L^2(\Omega)}^2.
\end{aligned}
\end{equation}
The augmentation terms vanish on the constraint manifold. Their first variations add a grad-div term to the velocity equation and a nonlinear radial term to the director equation. On the constraint manifold, the corresponding director Jacobian contribution is positive semidefinite in the normal direction; away from the manifold, an additional term proportional to \(|\mathbf n|^2-1\) is present and need not be positive. These contributions motivate the Schur-complement approximations developed below.

\subsection{Variational weak form of the augmented Lagrangian}
The augmented Lagrangian weak formulation is obtained by adding the penalty 
terms directly to the weak form of \eqref{eq:EL-saddle}. Recall that the 
original Lagrangian functional enforces the constraints 
weakly through the multipliers $p$ and $q$. The augmented Lagrangian 
functional further introduces quadratic penalty terms 
$\frac{\gamma_u}{2}\|\nabla\cdot\mathbf{u}\|^2$ and 
$\frac{\gamma_n}{2}\|\,|\mathbf{n}|^2-1\|^2$. 
At the level of the weak form, these penalty terms contribute 
$\gamma_u(\nabla\cdot\mathbf{u},\nabla\cdot\mathbf{v})$ to the momentum 
equation and 
$2\gamma_n(|\mathbf{n}|^2\mathbf{n},\mathbf{m})-2\gamma_n(\mathbf{n},\mathbf{m})$ 
to the director equation. With the homogeneous boundary conditions 
$\mathbf{u}|_{\partial\Omega}=0$ and $\partial_{\boldsymbol{\nu}}\mathbf{n}|_{\partial\Omega}=0$, 
the augmented weak form reads:

Find $\mathbf{u}\in V_u$, $\mathbf{n}\in V_n$, $p\in Q$ and $q\in\Lambda$ such that, for all test functions $\mathbf{v}\in V_u$, $\mathbf{m}\in V_n$, $r\in Q$ and $\zeta\in\Lambda$,

\begin{equation}\label{eq:aug-weak-u}
\begin{aligned}
&\left( \mathbf{u}_t, \mathbf{v} \right) + \left( \mathbf{u}\cdot\nabla \mathbf{u}, \mathbf{v} \right) - (p, \nabla\cdot \mathbf{v}) + \nu(\nabla \mathbf{u}, \nabla \mathbf{v})\\
&- K\left( (\nabla \mathbf{n})^\top\nabla \mathbf{n}, \nabla \mathbf{v} \right) + \gamma_u(\nabla\cdot \mathbf{u}, \nabla\cdot \mathbf{v}) = 0,
\end{aligned}
\end{equation}
\begin{equation}\label{eq:aug-weak-n}
\begin{aligned}
&\left( \mathbf{n}_t, \mathbf{m} \right) + \left( \mathbf{u}\cdot\nabla \mathbf{n}, \mathbf{m} \right) + \mu(\nabla \mathbf{n}, \nabla \mathbf{m}) + \mu\left(q \mathbf{n}, \mathbf{m} \right)\\
& + 2\gamma_n\left( |\mathbf{n}|^2 \mathbf{n}, \mathbf{m} \right) - 2\gamma_n\left( \mathbf{n}, \mathbf{m} \right)  = 0,
\end{aligned}
\end{equation}
\begin{equation}\label{eq:aug-weak-p}
(r, \nabla\cdot \mathbf{u}) = 0,
\end{equation}
\begin{equation}\label{eq:aug-weak-q}
(\zeta, |\mathbf{n}|^2 - 1) = 0,
\end{equation}
where the function spaces are defined as:
\begin{align*}
V_u&=H^1_0(\Omega)^d,\qquad V_n=H^1(\Omega)^d,\\
Q&=L^2_0(\Omega),\qquad \Lambda=L^2(\Omega).
\end{align*}
Here \(L^2_0(\Omega)\) fixes the additive constant in the pressure. The homogeneous Neumann boundary condition for the director is a natural boundary condition in the weak form and is therefore not imposed as an essential constraint on \(V_n\).

Compared with the unaugmented weak form of \eqref{eq:EL-saddle}, the formulation above adds the grad--div augmentation in the velocity 
equation and the nonlinear length-constraint augmentation in the director equation. The multiplier equations 
\eqref{eq:aug-weak-p}--\eqref{eq:aug-weak-q} are unchanged, so the constraints continue to be imposed weakly rather than replaced by a pure penalty approximation.

At the continuous level, the augmented and unaugmented constrained
problems are equivalent for sufficiently regular solutions satisfying the
constraints pointwise. Indeed, if
\(\nabla\cdot\mathbf{u}=0\) and \(|\mathbf{n}|=1\) almost everywhere, then
\[
\gamma_u(\nabla\cdot\mathbf{u},\nabla\cdot\mathbf{v})=0,
\qquad
2\gamma_n\bigl((|\mathbf{n}|^2-1)\mathbf{n},\mathbf{m}\bigr)=0
\]
for all admissible test functions \(\mathbf{v}\) and \(\mathbf{m}\).  
Hence the same quadruple also satisfies the augmented weak form
\eqref{eq:aug-weak-u}--\eqref{eq:aug-weak-q}.  Conversely, if a solution of the augmented weak form satisfies the two multiplier equations, 
then the augmentation terms vanish and the unaugmented saddle-point formulation is recovered.  
This elementary observation establishes equivalence at the continuous level. At the discrete level, parameter independence is obtained only when each augmentation is formed from the same finite-dimensional residual as its multiplier equation. The distinction is important here because a weak finite element constraint does not, by itself, make an unprojected pointwise penalty vanish; this issue is made explicit in Section~\ref{subsec:algebraic}.

\section{Finite element space and time discretizations}
\label{sec:Finite}

In this section, we describe the spatial discretization using the finite element method and the temporal discretization using the backward Euler method. Based on the augmented Lagrangian weak form derived in Section~\ref{sec:Augmented}, we present the fully discrete nonlinear system, which is then linearized by Newton's method. The resulting saddle-point system will serve as the foundation for the preconditioning techniques developed in Section~\ref{sec:Pre}.

\subsection{Finite element discretization}
The continuous function spaces \(V_u\),\(V_n\),\(Q\) and \(\Lambda\) are defined as in Section~\ref{sec:Augmented}.
We choose conforming finite element subspaces \(V_{u,h}\subset V_u\), \(V_{n,h}\subset V_n\), \(Q_h\subset Q\) and \(\Lambda_h\subset\Lambda\).  In this work we employ Taylor--Hood elements \(\mathbb{P}_2^d\)--\(\mathbb{P}_1\) for the velocity and pressure on simplicial meshes. For the director and the Lagrange multiplier we adopt the same pairing, namely $\mathbb{P}_2^d$ for $n_h$ and $\mathbb{P}_1$ for $q_h$. The pressure is taken in the mean-zero space in order to fix the additive constant. 
 
The unit-length constraint forms a saddle-point structure for the director and its multiplier. Linearizing the constraint about a discrete director state \(n_h\) gives the zero-order bilinear form
\begin{equation}
b_{n,h}(m_h,\zeta_h) := (2 n_h \cdot m_h,\zeta_h), \qquad \forall m_h\in V_{n,h},\; \zeta_h\in\Lambda_h.
\end{equation}

The discrete stability of this zero-order constraint is governed by the director--multiplier inf-sup condition
\begin{equation}
\inf_{\zeta_h\in\Lambda_h}
\sup_{m_h\in V_{n,h}}
\frac{(2n_h\cdot m_h,\zeta_h)}
{\|m_h\|_{n,h}\,\|\zeta_h\|_{L^2}}
\ge \beta_n ,
\label{eq:director_infsup}
\end{equation}
where \(\|\cdot\|_{n,h}\) denotes the mass-stiffness norm induced by the director block and \(\beta_n>0\) is independent of the mesh size. The standard sufficient assumptions are that the meshes are shape regular and quasi-uniform, and that the Newton iterates satisfy
\[
0<c_n\le |n_h(x)|\le C_n,\qquad \|n_h\|_{W^{1,\infty}(\Omega)}\le C_n,
\]
with constants independent of \(h\). Under these assumptions, a Fortin-type lifting can be constructed by projecting the continuous test direction \(\zeta_h n_h\) into the discrete director space. This gives a uniformly bounded test function \(m_h\) for which \(2(n_h\cdot m_h,\zeta_h)\) controls \(\|\zeta_h\|_{L^2}^2\). This is precisely the discrete tangent-space stability mechanism used in saddle-point finite element methods for harmonic maps and liquid-crystal flows; see, for example, the Ericksen--Leslie analysis of Badia et al.~\cite{Badia2011} and the related harmonic-map analyses in~\cite{Hu2009,Guti2017}.

The choice \(V_{n,h}=\mathbb{P}_2^d\) and \(\Lambda_h=\mathbb{P}_1\) is consistent with this stability mechanism: \((\mathbb{P}_1)^d\subset(\mathbb{P}_2)^d\), so the director test space contains the lower-order space used in established \(\mathbb{P}_1^d\)--\(\mathbb{P}_1\) analyses, as also explained in~\cite{xia2021}. We use this observation as the compatibility rationale for the present discretisation; a standalone proof of a uniform inf-sup constant for every possible Newton iterate is outside the scope of this paper, and the numerical tests below monitor the resulting constraint residuals directly.

\subsection{Algebraic matrix representation}
\label{subsec:algebraic}

Let \(\{\phi_j^u\}_{j=1}^{N_u}\), \(\{\phi_j^n\}_{j=1}^{N_n}\), \(\{\psi_j\}_{j=1}^{N_p}\) and \(\{\theta_j\}_{j=1}^{N_q}\) be bases of \(V_{u,h}\), \(V_{n,h}\), \(Q_h\) and \(\Lambda_h\) respectively. With this basis selection, the discrete expressions are derived as follows:
\[
u_h = \sum_{j=1}^{N_u} U_j\phi_j^u,\quad
n_h = \sum_{j=1}^{N_n} N_j\phi_j^n,\quad
p_h = \sum_{j=1}^{N_p} P_j\psi_j,\quad
q_h = \sum_{j=1}^{N_q} Q_j\theta_j,
\]
where \(U_j, N_j, P_j, Q_j\) are the time-dependent coefficients. Substituting these discrete functions into the augmented Lagrangian weak form \eqref{eq:aug-weak-u}--\eqref{eq:aug-weak-q} and choosing the test functions from the same discrete spaces, we obtain the semi-discrete problem: find \( (u_h, n_h, p_h, q_h) \in V_{u,h} \times V_{n,h} \times Q_h \times \Lambda_h \) such that for all test functions \( (v_h, m_h, r_h, \zeta_h) \) in the same spaces, the following equations hold:

\begin{equation}
\begin{aligned}
(u_{h,t},v_h) + (u_h\cdot\nabla u_h,v_h) - (p_h,\nabla\cdot v_h) + \nu(\nabla u_h,\nabla v_h) \\
- K\big((\nabla n_h)^\top\nabla n_h,\nabla v_h\big) + \gamma_u(\nabla\cdot u_h, \nabla\cdot v_h) = 0,
\end{aligned}
\end{equation}
\begin{equation}
\begin{aligned}
(n_{h,t},m_h) + (u_h\cdot\nabla n_h,m_h) + \mu(\nabla n_h,\nabla m_h) + \mu(q_hn_h,m_h)\\
+ 2\gamma_n(|n_h|^2 n_h,m_h) - 2\gamma_n(n_h,m_h)  = 0,
\end{aligned}
\end{equation}
\begin{equation}
(r_h,\nabla\cdot u_h) = 0,
\end{equation}
\begin{equation}
(\zeta_h,|n_h|^2-1) = 0.
\end{equation}

There is an analogous, and more consequential, distinction for the director constraint. The discrete multiplier equation
\[
(\zeta_h,|n_h|^2-1)=0\qquad\forall\zeta_h\in\Lambda_h
\]
annihilates only the \(\Lambda_h\)-moments of the nonlinear residual; it does not imply \(|n_h|^2-1=0\) pointwise. Consequently, the unprojected term
\[
2\gamma_n\bigl((|n_h|^2-1)n_h,m_h\bigr)
\]
need not vanish at an exactly solved discrete saddle point, and the discrete solution can depend on \(\gamma_n\). A finite-dimensional augmented Lagrangian that is exactly compatible with the multiplier equation would instead augment the algebraic residual, for example by
\[
\frac{\gamma_n}{2}\,\mathbf R(\mathbf n)^\top\mathbf M_q^{-1}\mathbf R(\mathbf n),
\]
whose first variation is \(\gamma_n\mathbf B_n^\top\mathbf M_q^{-1}\mathbf R(\mathbf n)\) and which vanishes whenever \(\mathbf R(\mathbf n)=0\).

To represent this system in matrix form, let $\mathbf{u}(t) \in \mathbb{R}^{N_u}$, $\mathbf{n}(t) \in \mathbb{R}^{N_n}$, $\mathbf{p}(t) \in \mathbb{R}^{N_p}$, and $\mathbf{q}(t) \in \mathbb{R}^{N_q}$ be the coefficient vectors. To keep the finite element formulas compact, discrete functions are written without boldface, whereas algebraic coefficient vectors and matrices are written in boldface. The mass matrices for velocity and director are denoted by $\mathbf{M}_u$ and $\mathbf{M}_n$ respectively, with entries
\[
(\mathbf{M}_u)_{ij}=(\phi_i^u,\phi_j^u),\qquad
(\mathbf{M}_n)_{ij}=(\phi_i^n,\phi_j^n).
\]

The stiffness matrices $\mathbf{K}_u$ and $\mathbf{K}_n$ are given by
\[
(\mathbf{K}_u)_{ij}=(\nabla\phi_i^u,\nabla\phi_j^u),\qquad
(\mathbf{K}_n)_{ij}=(\nabla\phi_i^n,\nabla\phi_j^n).
\]

The convection operators $\mathbf{C}_u(\mathbf{u})$ and $\mathbf{C}_n(\mathbf{u})$ are defined through their action:
\[
[\mathbf{C}_u(\mathbf{u})\mathbf{u}]_i = (u_h\cdot\nabla u_h,\phi_i^u),\qquad
[\mathbf{C}_n(\mathbf{u})\mathbf{n}]_i = (u_h\cdot\nabla n_h,\phi_i^n).
\]

The divergence matrix $\mathbf{B}$ is defined by \[
(\mathbf{B})_{ij}=(\psi_i,\nabla\cdot\phi_j^u).\]

The grad-div stabilization matrix $\mathbf{M}_{\text{div}}$ is defined by $(\mathbf{M}_{\text{div}})_{ij} = (\nabla\cdot\phi_i^u, \nabla\cdot\phi_j^u)$.

The pressure mass matrix $\mathbf{M}_p$ is defined by $(\mathbf{M}_p)_{ij}=(\psi_i,\psi_j)$. The mass matrix for the Lagrange multiplier is denoted by \(\mathbf{M}_{q}\) and has entries
$(\mathbf{M}_{q})_{ij} = (\theta_{i}, \theta_{j})$.

The nonlinear elastic coupling is represented by the residual operator
\[
[\mathbf{F}_{un}(\mathbf{n})]_i = -K\big((\nabla n_h)^\top\nabla n_h,\nabla\phi_i^u\big).
\]

For the director equation, the Lagrange multiplier coupling and the penalty nonlinearity are given by
\[[\mathbf{L}(\mathbf{n})\mathbf{q}]_i = (q_h n_h,\phi_i^n),\qquad
[\mathbf{H}_n(\mathbf{n})\mathbf{n}]_i = (|n_h|^2n_h,\phi_i^n).\]

Finally, the discrete constraint residual and the linearized constraint operator are
\[
[\mathbf{N}(\mathbf{n})\delta\mathbf{n}]_{i}=(n_h\cdot\delta n_h,\theta_i),\qquad
[\mathbf{R}(\mathbf{n})]_i = (|n_h|^2-1,\theta_i).
\]

With these definitions, the semi-discrete problem can be written compactly as
\begin{equation}
\label{eq:semi-discrete}
\begin{cases}
    &\mathbf{M}_u \dot{\mathbf{u}} + \mathbf{C}_u(\mathbf{u})\mathbf{u} + \nu \mathbf{K}_u \mathbf{u} - \mathbf{B}^\top \mathbf{p} + \gamma_u \mathbf{M}_{\text{div}} \mathbf{u} + \mathbf{F}_{un}(\mathbf{n}) = \mathbf{0},\\
&\mathbf{M}_n \dot{\mathbf{n}} + \mathbf{C}_n(\mathbf{u})\mathbf{n} + \mu \mathbf{K}_n \mathbf{n} + \mu \mathbf{L}(\mathbf{n})\mathbf{q} + 2\gamma_n (\mathbf{H}_n(\mathbf{n})\mathbf{n} - \mathbf{M}_n \mathbf{n}) = \mathbf{0},\\
&\mathbf{B} \mathbf{u} = \mathbf{0},\\
&\mathbf{R}(\mathbf{n}) = \mathbf{0}.
\end{cases}
\end{equation}

\subsection{Implicit Euler time integration}

Divide the time interval \([0,T]\) into \(N\) equal steps of size \(\Delta t = \frac{T}{N}\), and denote by \((u_h^n, n_h^n, p_h^n, q_h^n)\) the approximation at time \(t_n=n\Delta t\). Applying the backward Euler method to the semi-discrete system \eqref{eq:semi-discrete} yields the fully discrete nonlinear system for the new time level \((u_h, n_h, p_h, q_h):=(u_h^{n+1}, n_h^{n+1}, p_h^{n+1}, q_h^{n+1})\):
\begin{equation}
\begin{aligned}
\Big(\frac{u_h-u_h^n}{\Delta t},v_h\Big) + (u_h\cdot\nabla u_h,v_h) - (p_h,\nabla\cdot v_h) + \nu(\nabla u_h,\nabla v_h) \\
- K\big((\nabla n_h)^\top\nabla n_h,\nabla v_h\big) + \gamma_u(\nabla\cdot u_h, \nabla\cdot v_h) = 0, \\[4pt]
\end{aligned}
\end{equation}
\begin{equation}
\begin{aligned}
\Big(\frac{n_h-n_h^n}{\Delta t},m_h\Big) + (u_h\cdot\nabla n_h,m_h) + \mu(\nabla n_h,\nabla m_h) + \mu(q_h n_h,m_h)\\
+ 2\gamma_n(|n_h|^2 n_h,m_h) - 2\gamma_n(n_h,m_h)  = 0, \\[4pt]
\end{aligned}
\end{equation}
\begin{equation}
\begin{aligned}
(r_h,\nabla\cdot u_h) = 0, \\[4pt]
\end{aligned}
\end{equation}
\begin{equation}
\begin{aligned}
(\zeta_h,|n_h|^2-1) = 0,
\end{aligned}
\end{equation}

Combining the finite element discretization and matrix representations discussed above, we can derive a fully discrete formulation, expressed as the following system of nonlinear algebraic equations:    

\begin{equation}
\label{eq:fully_discrete}
\begin{cases}
\frac{1}{\Delta t}\mathbf{M}_u(\mathbf{u} - \mathbf{u}^n)
+ \mathbf{C}_u(\mathbf{u})\mathbf{u} + \nu\mathbf{K}_u\mathbf{u}
- \mathbf{B}^\top \mathbf{p}
+ \gamma_u\mathbf{M}_{\text{div}}\mathbf{u}
+ \mathbf{F}_{un}(\mathbf{n}) = \mathbf{0}, \\
\frac{1}{\Delta t}\mathbf{M}_n(\mathbf{n} - \mathbf{n}^n)
+ \mathbf{C}_n(\mathbf{u})\mathbf{n} + \mu\mathbf{K}_n\mathbf{n}
+ \mu\mathbf{L}(\mathbf{n})\mathbf{q} + 2\gamma_n\mathbf{H}_n(\mathbf{n})\mathbf{n}  - 2\gamma_n\mathbf{M}_n\mathbf{n} = \mathbf{0}, \\
\mathbf{B}\mathbf{u} = \mathbf{0}, \\
\mathbf{R}(\mathbf{n}) = \mathbf{0}.
\end{cases}
\end{equation}

\subsection{Newton linearization}\label{sec:newton}
The fully discrete system is nonlinear because of convection, elastic coupling, the multiplier--director product, and the director augmentation. We solve it with a damped Newton method and derive the exact Jacobian used at each iteration.

Denote the current state $(\mathbf{u}^k, \mathbf{n}^k, \mathbf{p}^k, \mathbf{q}^k)$, and the step increments $\delta\mathbf{u} = \mathbf{u}^{k+1} - \mathbf{u}^k$, $\delta\mathbf{n} = \mathbf{n}^{k+1} - \mathbf{n}^k$, $\delta\mathbf{p} = \mathbf{p}^{k+1} - \mathbf{p}^k$, $\delta\mathbf{q} = \mathbf{q}^{k+1} - \mathbf{q}^k$.

For the momentum equation, the linearized augmented system reads
\begin{equation}
\begin{split}
\left[
\frac{1}{\Delta t}\mathbf{M}_u+\nu \mathbf{K}_u
+\mathbf{C}_u(\mathbf{u}^k)
+\mathbf{C}_u'(\mathbf{u}^k)\mathbf{u}^k
+\gamma_u \mathbf{M}_{\text{div}}
\right]\delta\mathbf{u}
+\mathbf{A}_{un}\delta\mathbf{n}
-\mathbf{B}^\top\delta\mathbf{p}
=\mathbf{r}_u^k.
\end{split}
\end{equation}
where \(\mathbf r_u^k\) is the negative momentum residual. The
off-diagonal derivative \(\mathbf A_{un}\) is the matrix assembled from
\[
a_{un}^{k}(\delta n_h,v_h)
=
-K\Big(
(\nabla\delta n_h)^\top\nabla n_h^k
+
(\nabla n_h^k)^\top\nabla\delta n_h,
\nabla v_h
\Big),
\]
which is the exact derivative of the Ericksen elastic-stress term.

For the director equation, we now state the exact Jacobian in bilinear form. At the current Newton state, define
\begin{align}
a_{nu}^{k}(\delta u_h,m_h)
&=(\delta u_h\cdot\nabla n_h^k,m_h),\\
a_{nq}^{k}(\delta q_h,m_h)
&=\mu(\delta q_h\,n_h^k,m_h),
\end{align}
and
\begin{align}
a_{nn}^{k}(\delta n_h,m_h)
={}&
\Delta t^{-1}(\delta n_h,m_h)
+(u_h^k\cdot\nabla\delta n_h,m_h)
+\mu(\nabla\delta n_h,\nabla m_h)
+\mu(q_h^k\delta n_h,m_h)
\nonumber\\
&+2\gamma_n((|n_h^k|^2-1)\delta n_h,m_h)
+4\gamma_n(((n_h^k\cdot\delta n_h)n_h^k),m_h).
\label{eq:exact_director_jacobian}
\end{align}
The last two terms are the exact Fr\'echet derivative of
\[
G(n_h;m_h)=2\gamma_n((|n_h|^2-1)n_h,m_h).
\]
Thus the director Newton equation is
\[
a_{nu}^{k}(\delta u_h,m_h)
+a_{nn}^{k}(\delta n_h,m_h)
+a_{nq}^{k}(\delta q_h,m_h)
=(r_n^k,m_h).
\]
Therefore, we have
\begin{equation}
\begin{split}
&\mathbf{A}_{nu}\delta\mathbf{u}
+ \big[\frac{1}{\Delta t}\mathbf{M}_n+\mu\mathbf{K}_n + \mathbf{C}_n(\mathbf{u}^k) + \mu\frac{\partial\mathbf{L}}{\partial\mathbf{n}}(\mathbf{n}^k)\mathbf{q}^k
\\ &+2\gamma_n\mathbf{H}_n(\mathbf{n}^k)+4\gamma_n\mathbf{H}_n'(\mathbf{n}^k)\mathbf{n}^k - 2\gamma_n\mathbf{M}_n\big]\delta\mathbf{n}+ \mu\mathbf{L}(\mathbf{n}^k)\delta\mathbf{q} = \mathbf{r}_n^k
\end{split}
\end{equation}
where \(\mathbf{A}_{nu} = \frac{\partial\mathbf{C}_n}{\partial\mathbf{u}}(\mathbf{u}^k)\mathbf{n}^k\) and \(\mathbf{r}_n^k\) is the director residual. 
The terms involving \(\mathbf{H}_n(\mathbf{n}^k)\) and \(\mathbf{H}_n'(\mathbf{n}^k)\mathbf{n}^k\) arise from the linearization of the penalty term 
\(2\gamma_n(|\mathbf{n}|^2\mathbf{n}, \mathbf{m})\) with respect to \(\mathbf{n}\), yielding contributions \(2\gamma_n(|\mathbf{n}^k|^2\delta\mathbf{n}, \mathbf{m})\) and \(4\gamma_n((\mathbf{n}^k\cdot\delta\mathbf{n})\mathbf{n}^k, \mathbf{m})\), respectively.

For the two constraints, we obtain
\begin{equation}
\begin{aligned}
&\mathbf{B} \delta\mathbf{u} = \mathbf{r}_p^k, 
\\
&2\mathbf{N}(\mathbf{n}^k) \delta\mathbf{n} = \mathbf{r}_q^k,
\end{aligned}
\end{equation}
where every right-hand side \(\mathbf r_\bullet^k\) is the negative residual of the corresponding nonlinear equation evaluated at the current Newton state.

The resulting linearized system can be written as the following $4\times4$ block system
\begin{equation}
\label{eq:saddle_system}
\begin{bmatrix}
\mathbf{A}_{uu} & \mathbf{A}_{un} & -\mathbf{B}_u^\top & \mathbf{0} \\
\mathbf{A}_{nu} & \mathbf{A}_{nn} & \mathbf{0} & \mathbf{A}_{nq} \\
\mathbf{B}_u & \mathbf{0} & \mathbf{0} & \mathbf{0} \\
\mathbf{0} & \mathbf{B}_n & \mathbf{0} & \mathbf{0}
\end{bmatrix}
\begin{bmatrix}
\delta\mathbf{u} \\ \delta\mathbf{n} \\ \delta\mathbf{p} \\ \delta\mathbf{q}
\end{bmatrix}
=
\begin{bmatrix}
\mathbf{r}_u^k \\ \mathbf{r}_n^k \\ \mathbf{r}_p^k \\ \mathbf{r}_q^k
\end{bmatrix},
\end{equation}

where the blocks are defined as
\begin{equation}\label{eq:Newton}
\begin{aligned}
\mathbf{A}_{uu} &= \frac{1}{\Delta t}\mathbf{M}_u + \nu \mathbf{K}_u+ \gamma_u \mathbf{M}_{\text{div}} + \mathbf{C}_u(\mathbf{u}^k) + \mathbf{C}_u'(\mathbf{u}^k)\mathbf{u}^k, \\[4pt]
\mathbf{A}_{nn} &= \frac{1}{\Delta t}\mathbf{M}_n + \mu \mathbf{K}_n +\mathbf{C}_n(\mathbf{u}^k) + \mu\frac{\partial \mathbf{L}}{\partial \mathbf{n}}(\mathbf{n}^k)\mathbf{q}^k\\
&\quad + 2\gamma_n \mathbf{H}_n(\mathbf{n}^k) + 4\gamma_n \mathbf{H}_n'(\mathbf{n}^k)\mathbf{n}^k - 2\gamma_n \mathbf{M}_n , \\[4pt]
\mathbf{A}_{un} &= \frac{\partial \mathbf{F}_{un}}{\partial \mathbf{n}}(\mathbf{n}^k),\qquad
\mathbf{A}_{nu} = \frac{\partial \mathbf{C}_n}{\partial \mathbf{u}}(\mathbf{u}^k)\mathbf{n}^k, \\[4pt]
\mathbf{B}_u &= \mathbf{B},\qquad
\mathbf{B}_n = 2\mathbf{N}(\mathbf{n}^k),\qquad \mathbf{A}_{nq}=\mu\mathbf{L}(\mathbf{n}^k).
\end{aligned}
\end{equation}
Here \(\operatorname{Mat}(a)\) denotes the matrix assembled from the
corresponding bilinear form. With the adopted bases,
\(\mathbf A_{nq}=(\mu/2)\mathbf B_n^\top\).

\section{Preconditioner}
\label{sec:Pre}

We now develop a block preconditioner for the Newton saddle-point system obtained in Section~\ref{sec:Finite}.

Let the unknowns be grouped into two vectors
\[
x_1 = \begin{bmatrix} \delta\mathbf{u} \\ \delta\mathbf{n} \end{bmatrix}, \qquad
x_2 = \begin{bmatrix} \delta\mathbf{p} \\ \delta\mathbf{q}\end{bmatrix}.
\]
The saddle point system can be written as
\begin{equation}
\begin{bmatrix}
\mathcal{F} & \mathcal{G} \\
\mathcal{H} & 0
\end{bmatrix}
\begin{bmatrix}
x_1 \\ x_2
\end{bmatrix}
=
\begin{bmatrix}
b_1 \\ b_2
\end{bmatrix},
\end{equation}
where
\begin{equation}
\mathcal{F} = \begin{bmatrix}
A_{uu} & A_{un} \\[2pt]
A_{nu} & A_{nn}
\end{bmatrix},\quad
\mathcal{G} = \begin{bmatrix}
-B_u^{T} & 0 \\[2pt]  
0 & A_{nq}
\end{bmatrix},\quad
\mathcal{H} = \begin{bmatrix}
B_u & 0 \\[2pt]
0 & B_n
\end{bmatrix}.
\end{equation}
and \(B_n = 2\mathbf{N}(\mathbf{n}^k)\), \(A_{nq} = \mu\mathbf{L}(\mathbf{n}^k)\).  Notice that
\begin{equation}
A_{nq} = \frac{\mu}{2}\,B_n^{\top},
\end{equation}
a relation that will be used to simplify the multiplier Schur complement.

Equation \eqref{eq:saddle_system} is a standard saddle-point system with a zero \((2,2)\) block. 
If \(\mathcal{F}\) and the Schur complement \(\mathcal{S} = -\mathcal{H}\mathcal{F}^{-1}\mathcal{G}\) are invertible, 
the ideal upper-triangular preconditioner is
\begin{equation}
\mathcal{P}_{\text{ideal}} = \begin{bmatrix}
\mathcal{F} & \mathcal{G} \\
0 & \mathcal{S}
\end{bmatrix}.
\end{equation}
For right preconditioning,
\(\mathcal A\mathcal P_{\mathrm{ideal}}^{-1}=I+N\), where \(N\) is
strictly block lower triangular and satisfies \(N^2=0\). Hence the minimal
polynomial has degree at most two and GMRES converges in at most two
iterations in exact arithmetic~\cite{saad1986}.
Since we cannot apply exactly \(\mathcal{F}^{-1}\) or \(\mathcal{S}^{-1}\), we replace them by suitable approximations.

The block $\mathcal{F}$ couples velocity and director. In the parameter regimes considered here, the diagonal blocks capture the dominant time-stepping, diffusion and augmentation contributions. We therefore use the block-diagonal approximation
\begin{equation}
\tilde{\mathcal{F}} = \begin{bmatrix}
\mathbf{A}_{uu} & 0 \\
0 & \mathbf{A}_{nn}
\end{bmatrix},
\end{equation}

Using \(\widetilde{\mathcal{F}}^{-1}\) in place of \(\mathcal{F}^{-1}\) gives the block-diagonal Schur-complement approximation
\begin{equation}
\mathcal{S} \approx \widetilde{\mathcal{S}} =
\begin{bmatrix}
B_u A_{uu}^{-1} B_u^{T} & 0 \\[4pt]
0 & -B_n A_{nn}^{-1} A_{nq}
\end{bmatrix}.
\end{equation}

Therefore, the block-diagonal primal approximation yields a decoupled Schur \emph{model} with pressure component
\begin{equation}
S_p = B_u A_{uu}^{-1} B_u^{\top},
\end{equation}
and signed multiplier component
\begin{equation}
S_q = -B_n A_{nn}^{-1} A_{nq}
     = -\frac{\mu}{2}\,B_n A_{nn}^{-1} B_n^{\top}.
\end{equation}
For the multiplier solve and the spectral diagnostics, we work with the corresponding positive block
\(
\widehat S_q:=-S_q=(\mu/2)B_n A_{nn}^{-1}B_n^\top
\).
We now construct approximations for \(S_p^{-1}\) and \(\widehat S_q^{-1}\).

\subsection{Approximation of the pressure Schur complement}
The augmented velocity block admits the splitting
\begin{equation}
\mathbf{A}_{uu} = \widehat{A}_u + \gamma_u \mathbf{M}_{\text{div}},
\end{equation}
where $\widehat{A}_u = \frac{1}{\Delta t}\mathbf{M}_u + \nu \mathbf{K}_u + \mathbf{C}_u(\mathbf{u}^k) + \mathbf{C}_u'(\mathbf{u}^k)\mathbf{u}^k$ collects the unaugmented terms and $\mathbf{M}_p$ is the pressure mass matrix.  

As demonstrated by Heister and Rapin~\cite{2012Heister}, the exact grad-div matrix $\mathbf{M}_{\text{div}}$ can be decomposed into an algebraic $L^2$-projection part ($B_u^\top \mathbf{M}_p^{-1} B_u$) and a stabilizing fluctuation part. Because the fluctuation term vanishes asymptotically as the mesh size $h \to 0$, it is mathematically justified to approximate $\mathbf{M}_{\text{div}} \approx B_u^\top \mathbf{M}_p^{-1} B_u$ for the purpose of preconditioning. This yields the computationally tractable approximation of the velocity block:
\begin{equation}
\mathbf{A}_{uu} \approx \widehat{A}_u + \gamma_u B_u^\top \mathbf{M}_p^{-1} B_u.
\end{equation}

Applying the Sherman--Morrison--Woodbury formula~\cite{Bacuta2006} to this approximated block gives the identity
\begin{equation}
S_p^{-1} \approx \gamma_u \mathbf{M}_p^{-1} + S_{p,0}^{-1}, \qquad S_{p,0} = B_u \widehat{A}_u^{-1} B_u^\top.
\end{equation}

To motivate an approximation of \(S_{p,0}^{-1}\), we consider the symmetric
generalized-Stokes principal part of \(\widehat{A}_u\),
\begin{equation}
\widehat{A}_{u,0}
=\frac{1}{\Delta t}\mathbf{M}_u+\nu \mathbf{K}_u,\qquad S_{p,0}^{(s)}=B_u \widehat{A}_{u,0}^{-1} B_u^\top.
\end{equation}
For stable velocity--pressure pairs and on the discrete mean-zero pressure space,
the Cahouet--Chabard approximation~\cite{cahouet1988,elman2014} provides a
spectrally equivalent approximation of $\bigl(S_{p,0}^{(s)}\bigr)^{-1}$.
\begin{equation}
\bigl(S_{p,0}^{(s)}\bigr)^{-1}
\simeq
\nu \mathbf{M}_p^{-1} + \frac{1}{\Delta t}\mathbf{K}_p^{-1},
\label{eq:cc_equivalence}
\end{equation}
where $(\mathbf{K}_p)_{ij} = (\nabla \psi_i, \nabla \psi_j)$ is the pressure stiffness matrix.

We next include the augmented-Lagrangian contribution in the generalized-Stokes
principal part. Define
\begin{equation}
S_{p,\gamma}^{(s)}
=
B_u
\left(
\widehat{A}_{u,0}
+
\gamma_u B_u^\top \mathbf{M}_p^{-1}B_u
\right)^{-1}
B_u^\top.
\end{equation}
Applying the Sherman--Morrison--Woodbury formula to this principal part gives
\begin{equation}
\left(S_{p,\gamma}^{(s)}\right)^{-1}
=
\left(S_{p,0}^{(s)}\right)^{-1}
+
\gamma_u \mathbf{M}_p^{-1}.
\label{eq:smw_pressure_schur_stokes}
\end{equation}
Combining this with \eqref{eq:cc_equivalence} yields
\begin{equation}
\left(S_{p,\gamma}^{(s)}\right)^{-1}
\simeq
(\gamma_u+\nu) \mathbf{M}_p^{-1}
+
\frac{1}{\Delta t}\mathbf{K}_p^{-1}.
\label{eq:augmented_stokes_pressure_equivalence}
\end{equation}
We obtain the pressure
Schur-complement approximation
\begin{equation}
\widetilde{S}_p^{-1}
=
(\gamma_u+\nu) \mathbf{M}_p^{-1}
+
\frac{1}{\Delta t} \mathbf{K}_p^{-1}.
\label{eq:pressure_schur_approx}
\end{equation}
\begin{remark}
For the symmetric augmented generalized-Stokes principal part, the spectral equivalence in \eqref{eq:augmented_stokes_pressure_equivalence} follows directly from the Cahouet--Chabard estimate and the Sherman--Morrison--Woodbury identity. More precisely, suppose that the Taylor--Hood pair is inf-sup stable on the discrete mean-zero pressure space and that
\[
c_1\left(\nu\mathbf{M}_p^{-1}+\frac1{\Delta t}\mathbf{K}_p^{-1}\right)
\le
\left(S_{p,0}^{(s)}\right)^{-1}
\le
c_2\left(\nu\mathbf{M}_p^{-1}+\frac1{\Delta t}\mathbf{K}_p^{-1}\right)
\]
with constants independent of \(h\), \(\nu\), and \(\Delta t\). Adding the same positive term \(\gamma_u\mathbf{M}_p^{-1}\) to all three operators and using \eqref{eq:smw_pressure_schur_stokes} gives
\[
\min\{c_1,1\}
\left((\gamma_u+\nu)\mathbf{M}_p^{-1}+\frac1{\Delta t}\mathbf{K}_p^{-1}\right)
\le
\left(S_{p,\gamma}^{(s)}\right)^{-1}
\le
\max\{c_2,1\}
\left((\gamma_u+\nu)\mathbf{M}_p^{-1}+\frac1{\Delta t}\mathbf{K}_p^{-1}\right).
\]
Thus \eqref{eq:pressure_schur_approx} is spectrally equivalent to the inverse Schur complement of the augmented generalized-Stokes principal part. The nonsymmetric convection and Newton terms are retained in the velocity block but are not used in this Schur approximation; for the full nonsymmetric Newton operator the statement should be interpreted as a principal-part approximation rather than a uniform symmetric positive-definite spectral-equivalence theorem.
\end{remark}

For sufficiently large \(\gamma_u\), the augmented-Lagrangian contribution
dominates and the pressure Schur-complement inverse is well represented by
the pressure mass inverse,
$\widetilde{\mathbf S}_p^{-1}\approx \gamma_u\mathbf M_p^{-1}$. For fixed \(\gamma_u\) and \(\nu\), in the small time-step limit
\(\Delta t\to0\), the inverse-Laplacian contribution $\frac{1}{\Delta t}\mathbf K_p^{-1}$
becomes dominant, and the approximation behaves like an inverse Laplacian.

\subsection{Approximation of the director-multiplier Schur complement}
The multiplier constraint operator \(B_n\) is a zero-order operator, arising
from the pointwise linearization of the unit-length constraint:
\begin{equation}
B_n \delta\mathbf{n} = 2 \mathbf{n}^k \cdot \delta\mathbf{n}.
\end{equation}
This zero-order nature distinguishes the multiplier Schur complement from the pressure Schur complement associated with the divergence operator.

Recall from equation \eqref{eq:Newton} in Section~\ref{sec:newton} that the augmented director block
has the form
\begin{equation}
\mathbf{A}_{nn}
= \frac{1}{\Delta t}\mathbf{M}_n + \mu \mathbf{K}_n
+ \mathbf{C}_n(\mathbf{u}^k)
+ \mu\frac{\partial\mathbf{L}}{\partial\mathbf{n}}(\mathbf{n}^k)\mathbf{q}^k
- 2\gamma_n\mathbf{M}_n
+ 2\gamma_n\mathbf{H}_n(\mathbf{n}^k) + 4\gamma_n\mathbf{H}_n'(\mathbf{n}^k)\mathbf{n}^k.
\end{equation}

For the construction of a preconditioner, the penalty contributions are
separated into two parts according to their algebraic structure. The term
\(4\gamma_n\mathbf{H}_n'(\mathbf{n}^k)\mathbf{n}^k\) is positive semi-definite
and represents the dominant penalty effect. Following the projection argument
of Xia et al.~\cite{xia2021}, its action on the director increment can be
approximated as
\begin{equation}
\bigl[4\gamma_n\mathbf{H}_n'(\mathbf{n}^k)\mathbf{n}^k\bigr]\,\delta\mathbf{n}
\;\approx\; \gamma_n B_n^\top \mathbf{M}_q^{-1} B_n\,\delta\mathbf{n}.
\end{equation}

The remaining combined penalty contribution \(2\gamma_n\mathbf{H}_n(\mathbf{n}^k)-2\gamma_n\mathbf{M}_n\) involves the
factor \((|\mathbf{n}^k|^2-1)\) and is therefore small near convergence. Its
influence is absorbed into the augmentation parameter \(\gamma_n\) in the
preconditioner. This simplification retains the essential spectral properties
required for robust convergence~\cite{xia2021}.

We now consider the unaugmented multiplier Schur complement
\begin{equation}
S_{q,0}=B_n\widehat{\mathbf{A}}_n^{-1}B_n^\top ,
\end{equation}
where \(\widehat{\mathbf{A}}_n\) denotes the unaugmented part of
\(\mathbf{A}_{nn}\), i.e.\ all terms except the penalty contributions. The subsequent discussion primarily serves to identify the relevant scalings and should be interpreted as a heuristic argument applied to the symmetric
principal part of \(\widehat{\mathbf{A}}_n\). In particular, the convection and Newton linearization terms are omitted. At the continuous level, transport is skew-adjoint for a divergence-free velocity with compatible boundary conditions; the discrete operator used here is generally nonsymmetric because incompressibility is imposed only weakly. We therefore treat transport and the remaining Newton terms as perturbations of the principal reaction--diffusion block~\cite{Mardal2011}.

Assume temporarily that the principal part of the director block is dominated
by
\begin{equation}
\widehat{\mathbf{A}}_n \approx \frac1{\Delta t}\mathbf{M}_n + \mu \mathbf{K}_n.
\end{equation}
A frozen-coefficient analysis then suggests that the inverse of the
unaugmented multiplier Schur complement combines a temporal reaction scale and
an elliptic stiffness scale. In a formal dual-to-primal
interpretation~\cite{Mardal2011,Schoberl1999}, this leads to the heuristic
structure
\begin{equation}
S_{q,0}^{-1}
\approx
\frac1{\Delta t}\mathbf{M}_q^{-1}
+
\mu \mathbf{M}_q^{-1}\mathbf{K}_q \mathbf{M}_q^{-1},
\end{equation}
where
\begin{equation}
(\mathbf{K}_q)_{ij}=(\nabla\theta_j,\nabla\theta_i)
\end{equation}
is the stiffness matrix on the multiplier space. This expression is used
solely to indicate the possible stiffness scaling in the inverse of the
unaugmented Schur complement; it does not claim uniform spectral equivalence
for the full Newton--Jacobian block.

In the augmented Lagrangian preconditioner, this stiffness contribution is not
resolved explicitly. Using the projection approximation
\(4\gamma_n\mathbf{H}_n'(\mathbf{n}^k)\mathbf{n}^k
\approx \gamma_n B_n^\top \mathbf{M}_q^{-1}B_n\) and applying the
Sherman--Morrison--Woodbury formula to the augmented block
\(\widehat{\mathbf{A}}_n + \gamma_n B_n^\top \mathbf{M}_q^{-1}B_n\),
we obtain the following Sherman--Morrison--Woodbury decomposition for the positive Schur block associated with this approximating operator, denoted by \(\widehat{S}_{q,\gamma}\):
\begin{equation}
\widehat{S}_{q,\gamma}^{-1}
=
\frac{2}{\mu}
\left(
S_{q,0}^{-1}+\gamma_n\mathbf{M}_q^{-1}
\right),
\end{equation}
where \(S_{q,0}=B_n\widehat{\mathbf{A}}_n^{-1}B_n^\top\). The prefactor \(\mu/2\) and the sign convention follow from the relation
\(\mathbf{A}_{nq} = \mu\mathbf{L}(\mathbf{n}^k) = \frac{\mu}{2}B_n^\top\) adopted in Section~\ref{sec:newton}. If the opposite sign convention is used for the multiplier block,
the sign of the Schur-complement preconditioner should be adjusted accordingly.

In the parameter regime targeted by the augmented-Lagrangian preconditioner~\cite{Benzi2006,Farrell2019}, the explicit term \(\gamma_n\mathbf{M}_q^{-1}\)
is designed to dominate the unresolved portion of the multiplier
Schur-complement inverse. Meanwhile, the temporal reaction scale
\(1/\Delta t\) provides the natural mass-matrix contribution from the
unaugmented part. Combining these observations yields the following practical,
penalty-dominated approximation:
\begin{equation}
\widetilde{S}_q^{-1}
=\frac{2}{\mu}\left(\gamma_n+\frac{c}{\Delta t}\right)\mathbf{M}_q^{-1}.
\label{eq:final_Sq}
\end{equation}
Here \(c=\mathcal O(1)\) accounts for the temporal reaction scaling as well as
constants arising from nondimensionalisation and the choice of finite element
spaces; in practice, we take \(c=1\).

For large \(\gamma_n\), or in the small-time-step regime where the temporal
reaction contribution \(1/\Delta t\) is dominant, the approximation reduces to
a simple mass-matrix scaling. We stress that \eqref{eq:final_Sq} is a practical augmented-Lagrangian approximation to the inverse of the positive multiplier Schur block, rather than a uniform spectral-equivalence result for the full nonsymmetric Newton operator. The numerical experiments in Section~\ref{sec:Numerical} assess the
robustness of the proposed approximation with respect to the mesh size, the
time step, and the penalty parameters.

The preceding approximation can be supported by a short principal-part argument.  
Let \(\widehat{\mathbf{A}}_n\) be the symmetric positive definite reaction--diffusion part of the director block and 
assume that \(B_n\) has full row rank on the discrete multiplier space.  
Define
\[
T_0=B_n\widehat{\mathbf A}_n^{-1}B_n^\top,
\qquad
T_\gamma=B_n(\widehat{\mathbf A}_n+\gamma_n B_n^\top\mathbf M_q^{-1}B_n)^{-1}B_n^\top.
\]
The Sherman--Morrison--Woodbury identity gives
\[
T_\gamma^{-1}=T_0^{-1}+\gamma_n\mathbf M_q^{-1}.
\]
Thus the augmented inverse Schur block contains an explicit multiplier-mass contribution whose size is controlled by \(\gamma_n\).  
If, for the reaction-dominated part of \(\widehat{\mathbf A}_n\), the unaugmented inverse satisfies the mesh-independent scaling
\(T_0^{-1}\lesssim \Delta t^{-1}\mathbf M_q^{-1}\) on the multiplier space, then
\((\gamma_n+c\Delta t^{-1})\mathbf M_q^{-1}\) has the correct leading scaling for \(T_\gamma^{-1}\).  
This argument is intentionally restricted to the symmetric principal part; the nonsymmetric convection and Newton-coupling terms are treated by 
the outer FGMRES iteration and are tested numerically below.

\section{Numerical results}
\label{sec:Numerical}
\begin{figure}[htbp]
  \centering
  \resizebox{\textwidth}{!}{\begin{tikzpicture}[%
  every node/.style={draw=black, thick, anchor=west},
  grow via three points={one child at (0.0,-0.7) and
  two children at (0.0,-0.7) and (0.0,-1.4)},
  edge from parent path={(\tikzparentnode.210) |- (\tikzchildnode.west)}]
  
  \node {Newton solver with line search}
    child { node {Krylov solver (FGMRES)}
      child { node {Block preconditioner (Schur complement)}
        child { node {Momentum and director block ($u, n$)}
          child { node {Direct solver (LU / MUMPS)} }
        }
        child [missing] {} 
        child { node {Pressure and multiplier block ($p, \lambda$) }
          child { node {Approximate pressure Schur complement}
            child { node {Mass matrix: Jacobi preconditioner} }
            child { node {Stiffness matrix: CG solver + GAMG} }
          }
          child [missing] {} 
          child [missing] {} 
          child { node {Approximate multiplier Schur complement}
            child { node {Mass matrix: Jacobi preconditioner} }
          }
        }
      }
    };
\end{tikzpicture}}
  \caption{Schematic of the solution algorithm.}
  \label{fig:solver}
\end{figure}
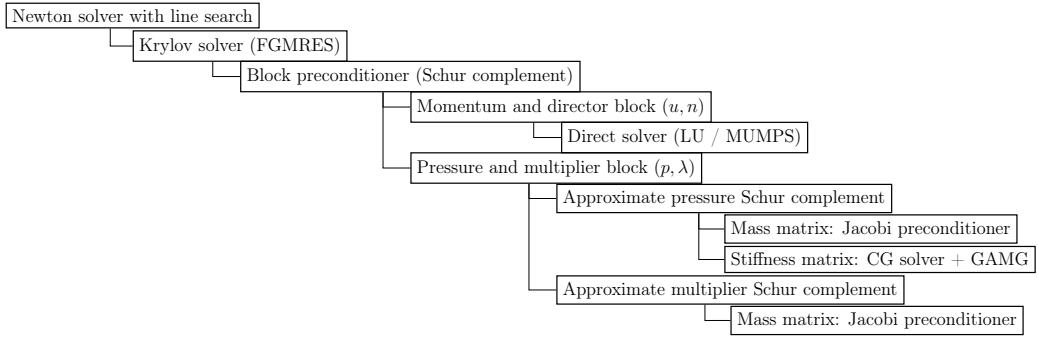

We assess the accuracy, structure preservation, and solver robustness of the fully discrete augmented Lagrangian formulation. The focus is on the Schur-complement-based preconditioner developed in Section~\ref{sec:Pre} for the Newton linearizations of the coupled Ericksen--Leslie saddle-point system. All numerical simulations are implemented using the open-source finite element framework Firedrake~\cite{firedrake2016} (version 2025.4.2) and PETSc (version 3.23.4), with time integration provided by the Irksome library~\cite{farrell2021irksome} (version 0.0.1). Computations were performed in serial (1 MPI rank) without threading on a virtual machine allocated with 28 GB of RAM, hosted on a workstation equipped with an Intel Core i7-13700H processor. A graphical representation of the entire algorithm is shown in Figure~\ref{fig:solver}.

Unless otherwise stated, the notation \(\gamma\) in this section means \(\gamma_u=\gamma_n=\gamma\). 
The numerical tests are organized as follows. A manufactured-solution test first verifies the convergence of the $\mathbb{P}_2$--$\mathbb{P}_1$ finite element discretization and the backward Euler time integrator. A smooth benchmark problem is then used to examine whether the augmented formulation is consistent with the expected energy-dissipation behavior. We next investigate the algebraic properties of the proposed block preconditioner. Subsequent tests evaluate the robustness of the global solver against variations in mesh size, time step, and fluid viscosity. Finally, computations on a disk and for a two-defect configuration illustrate the applicability of the method to curved geometries and director fields with large gradients.

In the reported implementation, the global nonlinear systems are solved using a Newton line-search method with a baseline relative tolerance of $10^{-8}$ and a baseline maximum of 50 iterations. The linearized saddle-point systems are solved using an outer FGMRES Krylov method with a baseline relative tolerance of $10^{-7}$. The coupled primal block \(\widetilde{\mathcal F}\), containing the velocity and director unknowns, is inverted using a sparse direct LU factorization with MUMPS. In the baseline benchmark and application computations, the pressure and multiplier mass inverses are approximated by single applications of Jacobi scaling, while the pressure-stiffness inverse is approximated using CG preconditioned with GAMG. For the manufactured-solution, algorithmic diagnostic, and robustness studies, the nonlinear and linear tolerances, as well as the auxiliary inner solves, are adjusted as needed to avoid contamination by algebraic errors, expose parameter dependence, or ensure accurate algebraic diagnostics. Thus the iteration studies below assess the quality of the outer block factorization and Schur approximations. Because the primal blocks are solved directly, they do not by themselves establish optimal-complexity or parallel scalability of the complete solver.

\subsection{Convergence analysis via the method of manufactured solutions}
We first verify the spatial and temporal accuracy of the fully discrete scheme using the method of manufactured solutions. The test is posed on the unit square $\Omega=(0,1)^2$. The manufactured solution is constructed 
so that the incompressibility constraint $\nabla\cdot\mathbf{u}=0$ and the unit-length constraint $|\mathbf{n}|=1$ hold pointwise.  The velocity is defined by a stream function $\psi$, the director by a rotation angle $\theta$, and smooth functions are prescribed for the pressure and the Lagrange multiplier:
\begin{align}
\psi(x,y,t) &= \sin^2(\pi x)\sin^2(\pi y)\cos(t), \qquad
\mathbf{u}_{\text{ex}} = \bigl(\partial_y\psi,\;-\partial_x\psi\bigr)^{\!\top},\\[2mm]
\theta(x,y,t) &= 0.4\cos(2\pi x)\cos(2\pi y)\cos(t), \qquad
\mathbf{n}_{\text{ex}} = (\cos\theta,\;\sin\theta)^{\!\top},\\[2mm]
p_{\text{ex}}(x,y,t) &= \sin(2\pi x)\cos(2\pi y)\sin(t+0.3),\\[2mm]
q_{\text{ex}}(x,y,t) &= 0.2\sin(\pi x)\sin(\pi y)\cos(t+0.1).
\end{align}
The corresponding forcing terms $\mathbf{f}_u$ and $\mathbf{f}_n$ are obtained by substituting the exact solution into the augmented weak form \eqref{eq:aug-weak-u}--\eqref{eq:aug-weak-q}. 
Physical parameters are fixed as $\nu=1.0$, $K=1.0$, $\mu=1.0$; the augmentation parameters are set to $\gamma_u=10.0$ and $\gamma_n=10.0$. 
Spatial discretization employs $\mathbb{P}_2$ elements for $\mathbf{u}$ and $\mathbf{n}$, and $\mathbb{P}_1$ elements for $p$ and $q$. 
Temporal discretization is performed with the backward Euler method. For these verification tests, the auxiliary mass and pressure-stiffness problems are solved using sparse direct LU factorizations with MUMPS, so that algebraic errors remain negligible relative to the measured discretization errors.

\subsubsection{Spatial convergence}
To reduce temporal contamination in the spatial study, we take $\Delta t=10^{-3}$ and advance only two time steps. The mesh is uniformly refined from $8\times8$ to $128\times128$ cells. Table~\ref{tab:mms_space} reports the $L^2$ errors for $\mathbf{u}$, $\mathbf{n}$, $p$, and $q$, together with the constraint residuals 
$\|\nabla\cdot\mathbf{u}\|_{L^2}$ and $\||\mathbf{n}|^2-1\|_{L^2}$.

\begin{table}[htbp]
  \centering
  \caption{Spatial convergence and constraint residuals for the manufactured-solution test on uniform meshes with $\mathbb{P}_2\!-\!\mathbb{P}_1$ elements and $\Delta t=10^{-3}$.} 
  \label{tab:mms_space}
  \begin{tabular}{@{}ccccccc@{}} 
\toprule
$N$ & $\|\mathbf{u}-\mathbf{u}_{\rm ex}\|_{L^2}$ & rate & $\|\mathbf{n}-\mathbf{n}_{\rm ex}\|_{L^2}$ & rate & $\|p-p_{\rm ex}\|_{L^2}$ & rate \\
\midrule
  8 & 1.37e-02 & -- & 1.70e-03 & -- & 3.05e-01 & -- \\
 16 & 1.54e-03 & 3.16 & 2.30e-04 & 2.89 & 2.74e-02 & 3.48 \\
 32 & 1.76e-04 & 3.13 & 2.96e-05 & 2.95 & 2.34e-03 & 3.55 \\
 64 & 2.12e-05 & 3.05 & 3.74e-06 & 2.98 & 2.26e-04 & 3.38 \\
128 & 2.63e-06 & 3.01 & 4.70e-07 & 2.99 & 3.40e-05 & 2.73 \\
\midrule\midrule 
$N$ & $\|q-q_{\rm ex}\|_{L^2}$ & rate & $\|\nabla\!\cdot\!\mathbf{u}\|_{L^2}$ & rate & $\||\mathbf{n}|^2-1\|_{L^2}$ & rate \\
\midrule
  8 & 1.57e-01 & -- & 3.38e-01 & -- & 9.00e-04 & -- \\
 16 & 2.01e-02 & 2.97 & 9.98e-02 & 1.76 & 1.35e-04 & 2.74 \\
 32 & 3.91e-03 & 2.36 & 2.67e-02 & 1.90 & 1.87e-05 & 2.85 \\
 64 & 1.16e-03 & 1.76 & 6.82e-03 & 1.97 & 2.41e-06 & 2.96 \\
128 & 1.13e-04 & 3.36 & 1.71e-03 & 2.00 & 3.04e-07 & 2.99 \\
\bottomrule
\end{tabular}
\end{table}

The velocity and director errors exhibit approximately third-order $L^2$ decay over the reported meshes, consistent with the expected behaviour of 
quadratic finite elements.  The pressure and multiplier errors also decrease under mesh refinement and attain at least the expected accuracy.  
In several refinement levels, their observed rates are locally higher than the nominal order associated with $\mathbb{P}_1$ approximations; this behaviour is discussed in the following remark.  
The two constraint residuals decrease systematically: the divergence residual is approximately second order, while the director-length residual is close to third order.  These results confirm the consistency of the augmented saddle-point discretisation for both the primary variables and the constraint variables, and show that the discrete incompressibility and unit-length constraints are increasingly well enforced under mesh refinement.

\begin{remark}
\label{rem:observed_high_order}
The pressure and multiplier errors in Table~\ref{tab:mms_space} show locally higher rates than the nominal $L^2$ accuracy expected for $\mathbb{P}_1$ approximations. We report this only as a numerical observation for the present manufactured-solution test, and do not interpret it as a general superconvergence result. The exact solution is smooth, the meshes are uniformly refined, and the computation is carried out over a very short time interval; under these favourable conditions, error cancellation and pre-asymptotic effects may lead to rates above the nominal order.
\end{remark}

\subsubsection{Temporal convergence}

We next verify the temporal accuracy of the backward Euler discretization on a fixed $64\times64$ mesh.  The time step is successively refined from $\Delta t=0.1$ to $\Delta t=6.25\times10^{-3}$, and the final time is $t_{\mathrm{end}}=0.1$.  The forcing terms are generated from the continuous time derivative of the manufactured solution.  Table~\ref{tab:mms_time} lists the $L^2$ errors and the corresponding convergence rates.
\begin{table}[htbp]
  \centering
  \caption{Temporal convergence on a fixed $64\times64$ mesh,
    $t_{\mathrm{end}}=0.1$.}
  \label{tab:mms_time}
  \footnotesize
  \begin{tabular}{@{}ccccccccc@{}}
    \toprule
    $\Delta t$ & $\|\mathbf{u}-\mathbf{u}_{\rm ex}\|_{L^2}$ & rate & $\|\mathbf{n}-\mathbf{n}_{\rm ex}\|_{L^2}$ & rate & $\|p-p_{\rm ex}\|_{L^2}$ & rate & $\|q-q_{\rm ex}\|_{L^2}$ & rate \\
\midrule
0.10000 & 1.53e-03 & -- & 1.09e-04 & -- & 5.39e-03 & -- & 3.65e-03 & -- \\
0.05000 & 8.46e-04 & 0.86 & 5.92e-05 & 0.89 & 2.95e-03 & 0.87 & 2.17e-03 & 0.75 \\
0.02500 & 4.44e-04 & 0.93 & 3.07e-05 & 0.95 & 1.54e-03 & 0.94 & 1.45e-03 & 0.58 \\
0.01250 & 2.30e-04 & 0.95 & 1.59e-05 & 0.95 & 7.97e-04 & 0.95 & 1.22e-03 & 0.25 \\
0.00625 & 1.20e-04 & 0.94 & 8.71e-06 & 0.87 & 4.40e-04 & 0.86 & 1.20e-03 & 0.03 \\
    \bottomrule
  \end{tabular}
\end{table}

The velocity, director, and pressure errors exhibit approximately first-order convergence, consistent with backward Euler time discretization.  The multiplier error decreases for the larger time 
steps, but its convergence rate degrades as $\Delta t$ becomes small.  
This behavior is expected on a fixed spatial mesh: $q$ is an algebraic 
constraint variable that does not appear under a time derivative, and 
once the temporal error in the director field is sufficiently reduced, 
the error in $q$ becomes dominated by the spatial discretization error.

\subsection{Energy dissipation}

To examine whether the computed solution exhibits the expected energy-dissipation trend of the Ericksen--Leslie model, we consider the smooth test problem described by Badia et al.~\cite{Badia2011}. 
The computational domain is the square $\Omega=[-1,1]^2$ with homogeneous Dirichlet boundary conditions for the velocity and homogeneous Neumann conditions for the director. The initial conditions and physical parameters are adopted from~\cite{Badia2011}:
\begin{equation}
\mathbf{u}_0 = \mathbf{0},\qquad 
\mathbf{n}_0 = (\sin a,\;\cos a)^{\!\top},\quad 
a = 2\pi\bigl(\cos x - \sin y\bigr),
\end{equation}
and the physical constants are set to $\nu=0.1$, $K=1.0$, $\mu=1.0$. The augmentation parameters are $\gamma_u=10.0$, $\gamma_n=10.0$, and the time step is $\Delta t = 2.5\times10^{-4}$; the simulation is run until $T=0.5$. Spatial discretization uses $\mathbb{P}_2$ elements for $\mathbf{u}$ and $\mathbf{n}$, and $\mathbb{P}_1$ elements for $p$ and $q$ on a uniform mesh of $50\times50$ cells (giving $h=1/25$).
The nonlinear system at each time step is solved with Newton's method. The linearized saddle-point system is solved using FGMRES preconditioned by the augmented Lagrangian Schur complement approximations derived in Section~\ref{sec:Pre}.

\begin{figure}[htbp]
\centering
\includegraphics[width=0.9\linewidth]{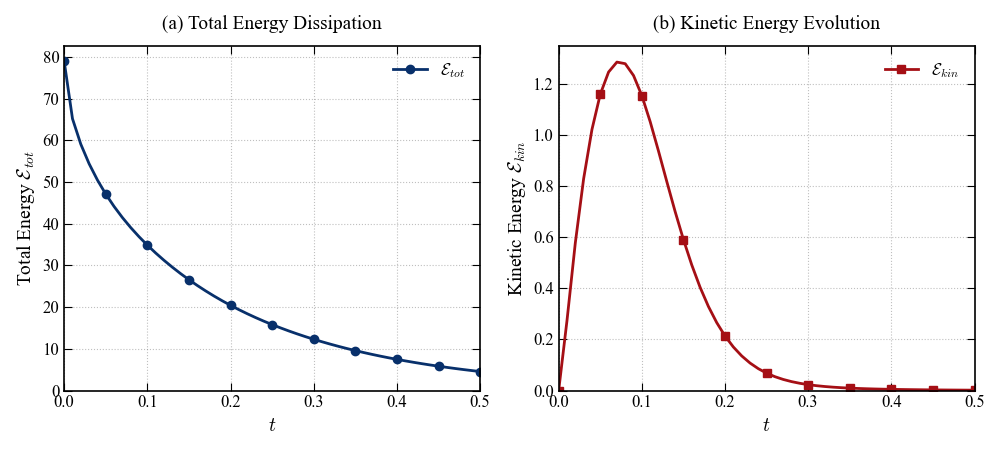}
\caption{Total energy and kinetic energy as functions of time for the smooth test.}
\label{fig:energy_badia}
\end{figure}

Figure~\ref{fig:energy_badia} shows the computed total energy $\mathcal{E}_{tol}(t)=\tfrac12\|u_h(t)\|^2+\tfrac K2\|\nabla n_h(t)\|^2$ and 
the kinetic energy $\mathcal{E}_{kin}(t)$. The total energy decreases monotonically at the sampled times, while the kinetic energy first rises 
and then decays. This behavior is consistent with the benchmark in~\cite{Badia2011} and with the continuous energy dissipation law.

\subsection{Mechanisms and spectral analysis of the block preconditioner}
\label{sec:Mechanisms}

In this section, we investigate the underlying mechanisms of the proposed augmented block preconditioner. Specifically, we compare simultaneous and single-constraint augmentations, examine the dominant influence of each penalty parameter on solver performance, and report algebraic diagnostics for the Schur models. 

Unless otherwise specified, all diagnostic tests in this section are performed on the unit square $\Omega=(0,1)^2$ with the smooth initial condition
\begin{equation}
\mathbf{u}_0 = \mathbf{0},\qquad 
\mathbf{n}_0 = (\sin a,\;\cos a)^{\!\top},\quad 
a = 2\pi(\cos x - \sin y),
\end{equation}
and with physical parameters $K = 1.0$, $\mu = 1.0$, $\Delta t = 10^{-3}$, 
$\nu = 0.1$, and $T = 0.02$. 
The parameter diagnostics are conducted on a $32\times32$ baseline mesh, and the spectral analysis uses a coarser $8\times8$ mesh.
The nonlinear systems are solved by Newton's method with a relative residual tolerance of $10^{-6}$. 
The linearized saddle-point systems are solved using preconditioned FGMRES with a relative tolerance of $10^{-4}$, where the preconditioner is based on the augmented Lagrangian Schur complement approximations presented in Section~\ref{sec:Pre}. The inner Conjugate Gradient solver for the pressure-stiffness block uses a relative tolerance of $10^{-10}$ to ensure that the algebraic error does not contaminate the parameter-dependence observations. In the augmentation comparison and parameter-sweep tests of this section, the pressure mass matrix and the multiplier mass matrix are treated with direct LU factorizations, rather than the Jacobi approximations used in the baseline configuration, so that the reported iteration counts are free of inner approximation errors.

\subsubsection{Comparison of augmentation strategies}

Four algorithmic variants are compared: an unaugmented baseline setting $\gamma_u=\gamma_n=0$, a pressure-only augmentation setting $\gamma_u=1000$ and $\gamma_n=0$, a multiplier-only augmentation setting $\gamma_u=0$ and $\gamma_n=1000$, and the proposed simultaneous augmentation setting both parameters to $1000$.

Table~\ref{tab:aug_comparison} reports the average Newton iterations per time step, the average FGMRES iterations per Newton step, the mean CPU time, and the continuous $L^2$ norms of the constraint residuals. The unaugmented baseline requires the highest number of Krylov iterations. Applying either the pressure-only or the multiplier-only augmentation independently results in a partial reduction in the average FGMRES iterations. This indicates that both the incompressibility constraint and the unit-length constraint contribute to the ill-conditioning of the linearized saddle-point system. When both augmentations are applied simultaneously, the Krylov iteration count drops to 5.35, resulting in the minimum overall computational cost.

Regarding constraint preservation, the director-length error remains essentially unchanged across all four configurations. 
Since the unit-length constraint is already enforced by the Lagrange multiplier up to the spatial truncation limit, the 
$\gamma_n$-augmentation functions mainly as an algebraic preconditioner rather than a geometric penalty. In contrast, the pressure augmentation 
actively decreases the divergence error. Consequently, the simultaneous strategy provides the dual advantage of maximizing solver efficiency while 
enhancing discrete mass conservation.

\begin{table}[htbp]
  \centering
  \caption{Comparison of augmentation strategies on a $32 \times 32$ mesh with $\Delta t=10^{-3}$ and $\nu=0.1$.}
  \label{tab:aug_comparison}
  \footnotesize
  \begin{tabular}{@{}lccccccc@{}}
    \toprule
    Variant & $\gamma_u$ & $\gamma_n$ & Newton it & KSP it & CPU (s) & $\|\nabla\cdot\mathbf{u}\|_{L^2}$ & $\||\mathbf{n}|^2-1\|_{L^2}$ \\
    \midrule
    No Augment    & 0 & 0 & 2.45 & 19.76 & 12.45 & 3.53e-02 & 3.00e-05 \\
    Pressure Only & 1000 & 0 & 2.55 & 10.20 & 8.63 & 2.69e-03 & 3.00e-05 \\
    Director Only & 0 & 1000 & 2.80 & 11.29 & 9.77 & 3.56e-02 & 2.99e-05 \\
    Simultaneous  & 1000 & 1000 & 2.85 & 5.35 & 7.70 & 2.69e-03 & 2.99e-05 \\
    \bottomrule
  \end{tabular}
\end{table}

\subsubsection{Independent effects of the augmentation parameters}

To further understand the mechanisms of the block preconditioner, we investigate how the grad--div parameter $\gamma_u$ and the director-length parameter $\gamma_n$ independently affect the solver convergence and constraint preservation. Specifically, one parameter is varied over four orders of magnitude while the other is kept fixed at a constant value of $100.0$. The numerical results are summarized in Table~\ref{tab:independent_gamma}.

Table~\ref{tab:independent_gamma} shows that increasing either augmentation parameter reduces the outer FGMRES count, with the stronger 
effect associated with the corresponding Schur approximation. The director-length diagnostic changes only marginally. 
These data support distinct dominant roles for the two
parameters over this test range. 
These tests confirm the independence of the two augmentations: each parameter improves the conditioning of its corresponding Schur block 
without interfering with the other constraint. This independence firmly justifies the block-diagonal design of the proposed preconditioner.

\begin{table}[htbp]
  \centering
  \caption{Independent effects of the augmentation parameters $\gamma_u$ and $\gamma_n$ on a $32\times32$ mesh with $\Delta t=10^{-3}$ and $\nu=0.1$.}
  \label{tab:independent_gamma}
  \begin{tabular}{@{}cccccc@{}}
    \toprule
    $\gamma_u$ & $\gamma_n$ & Newton it & KSP it & $\|\nabla\cdot\mathbf{u}_h\|_{L^2}$ & $\||\mathbf{n}_h|^2-1\|_{L^2}$ \\
    \midrule
    1.0    & 100.0 & 2.65 & 8.66 & 2.33e-02 & 3.00e-05 \\
    10.0   & 100.0 & 2.55 & 6.73 & 1.85e-02 & 3.00e-05 \\
    100.0  & 100.0 & 2.55 & 6.12 & 9.42e-03 & 3.00e-05 \\
    1000.0 & 100.0 & 3.00 & 5.92 & 2.69e-03 & 3.00e-05 \\
    \midrule
    100.0 & 1.0    & 2.75 & 7.07 & 9.42e-03 & 3.00e-05 \\
    100.0 & 10.0   & 2.75 & 7.04 & 9.42e-03 & 3.00e-05 \\
    100.0 & 100.0  & 2.55 & 6.12 & 9.42e-03 & 3.00e-05 \\
    100.0 & 1000.0 & 2.45 & 4.06 & 9.42e-03 & 2.99e-05 \\
    \bottomrule
  \end{tabular}
\end{table}

\subsubsection{Spectral analysis of the preconditioned Schur complements}

We next report a small-matrix spectral diagnostic for the diagonal Schur model induced by \(\widetilde{\mathcal F}\). On an $8\times8$ mesh, the Jacobian is assembled about a frozen, nonphysical state with constant velocity $\mathbf u=(0.1,-0.1)^\top$ so that the transport blocks are nonzero. We form \(S_p=B_uA_{uu}^{-1}B_u^\top\) and the positive multiplier block \(\widehat S_q=(\mu/2)B_n A_{nn}^{-1}B_n^\top\). These are exact blocks of the block-diagonal Schur model, not the diagonal blocks of the full coupled Schur complement \(-\mathcal H\mathcal F^{-1}\mathcal G\). For the reported matrices, the computed spectra were real and positive; Table~\ref{tab:spectral_diagnostics} lists their extreme eigenvalues and the eigenvalue-spread ratio \(\rho_{\rm eig}=\lambda_{\max}/\lambda_{\min}\). For a general nonsymmetric Newton state, singular values or field-of-values information would be more appropriate than this ratio.

\begin{table}[htbp]
  \centering
  \caption{Extreme real eigenvalues and eigenvalue-spread ratios \(\rho_{\rm eig}=\lambda_{\max}/\lambda_{\min}\) for the preconditioned diagonal Schur-model blocks on an \(8\times8\) mesh.}
  \label{tab:spectral_diagnostics}
  \begin{tabular}{@{}ccccccc@{}}
    \toprule
    & \multicolumn{3}{c}{Pressure block $\widetilde{S}_p^{-1} S_p$} & \multicolumn{3}{c}{Multiplier block $\widetilde{S}_q^{-1} \widehat S_q$} \\
    \cmidrule(lr){2-4} \cmidrule(l){5-7}
    $\gamma$ & $\lambda_{\min}$ & $\lambda_{\max}$ & $\rho_{\rm eig}$ & $\lambda_{\min}$ & $\lambda_{\max}$ & $\rho_{\rm eig}$ \\
    \midrule
    1.0    & 0.24 & 0.97 & 4.07 & 1.56 & 3.88 & 2.50 \\
    10.0   & 0.22 & 0.98 & 4.47 & 1.55 & 3.79 & 2.45 \\
    100.0  & 0.22 & 1.00 & 4.53 & 1.48 & 3.08 & 2.09 \\
    1000.0 & 0.22 & 1.00 & 4.47 & 1.20 & 1.59 & 1.32 \\
    \bottomrule
  \end{tabular}
\end{table}

Table~\ref{tab:spectral_diagnostics} reports the spectral bounds and condition numbers for varying values of the augmentation parameter. 
For the multiplier block, the condition number decreases monotonically from 2.50 to 1.32 as $\gamma$ grows.  This behavior is consistent with the structure of the Schur complement approximation~\eqref{eq:final_Sq}: for large $\gamma$, the penalty term $\gamma\,\mathbf{M}_q$ dominates the exact multiplier Schur complement $S_q$, making it spectrally close to a scalar multiple of the mass matrix.  Consequently, the preconditioner $\widetilde{S}_q^{-1}$ becomes an increasingly accurate approximation, leading to tightly clustered eigenvalues and a reduced condition number.

The condition number of the preconditioned pressure block remains uniformly bounded at approximately $4.5$, independent of $\gamma$.  In the exact system, the inclusion of the grad-div penalty introduces algebraic stiffness into $S_p$.  The $\gamma$-independence observed in Table~\ref{tab:spectral_diagnostics} demonstrates that this stiffness is effectively neutralized by the scaled mass matrix component $(\gamma_u + \nu)\mathbf{M}_p^{-1}$ in the Schur complement approximation~\eqref{eq:pressure_schur_approx}.  The eigenvalues remain bounded within $[0.22, 1.00]$, ensuring that the spectral equivalence does not deteriorate for large augmentation parameters.

These spectral diagnostics confirm that the proposed block preconditioner correctly resolves the different scalings of both constraints.  The asymptotic exactness of the multiplier approximation and the uniform boundedness of the preconditioned pressure block prevent local spectral degradation.

\subsection{Robustness of the augmented Lagrangian preconditioner}
\label{sec:robustness}

We now assess the outer iteration counts with respect to mesh size, the augmentation parameter $\gamma$, the time step, and the fluid viscosity $\nu$. 
Unless stated otherwise, the setup and nonlinear/linear tolerances are those of Section~\ref{sec:Mechanisms}. 
The reported study concerns outer Newton--FGMRES robustness.
The solution algorithm, including the inner solvers, is implemented as 
depicted in Figure~\ref{fig:solver}; all inner solves are performed to sufficiently high accuracy so that their contribution to the overall error 
is negligible.

\subsubsection{Mesh refinement and augmentation parameter}

We first examine the dependence of the solver on the mesh size and the augmentation parameter $\gamma$.  
Table~\ref{tab:mesh_gamma} reports the performance for four tested mesh resolutions and four values of $\gamma$.  
The viscosity is kept constant at $\nu = 0.1$.

For fixed $\gamma$, the average FGMRES count remains bounded with only mild growth over the tested mesh range. For a fixed mesh, 
increasing $\gamma$ reduces the count from approximately $9$--$11$ at $\gamma=1$ to approximately $3$--$4$ at $\gamma=1000$, consistent with 
augmentation-dominated Schur scaling.  
The reported divergence error $\|\nabla\cdot\mathbf{u}_h\|_{L^2}$ generally decreases with both large $\gamma$ and mesh refinement, 
reflecting the enhanced enforcement of incompressibility via the penalisation. 
The director-length diagonastic is governed mainly by mesh refinement in this test, 
because the length constraint is already enforced through the multiplier equation.

\begin{table}[htbp]
\centering
\caption{Robustness with respect to mesh refinement and $\gamma$. }
\label{tab:mesh_gamma}
\begin{tabular}{ccccccc}
\toprule
Mesh & DoFs & $\gamma$ & Newton it & KSP it & $\|\nabla\cdot\mathbf{u}\|_{L^2}$ & $\||\mathbf{n}|^2-1\|_{L^2}$ \\
\midrule
$8\times8$   & 1318     & 1     & 2.95 & 9.97 & 3.02e-01 & 1.71e-03 \\
             &          & 10    & 2.55 & 8.22 & 1.26e-01 & 1.70e-03 \\
             &          & 100   & 2.35 & 5.53 & 2.73e-02 & 1.67e-03 \\
             &          & 1000  & 2.40 & 3.75 & 3.24e-03 & 1.58e-03 \\
\addlinespace
$32\times32$  & 19078   & 1     & 2.75 & 9.42 & 2.33e-02 & 3.00e-05 \\
             &          & 10    & 2.80 & 7.29 & 1.85e-02 & 3.00e-05 \\
             &          & 100   & 2.55 & 5.96 & 9.42e-03 & 3.00e-05 \\
             &          & 1000  & 2.55 & 3.25 & 2.69e-03 & 2.99e-05 \\
\addlinespace
$64\times64$  & 75014   & 1     & 2.80 & 10.09 & 6.08e-03 & 3.80e-06 \\
             &          & 10    & 2.85 & 7.70 & 5.50e-03 & 3.80e-06 \\
             &          & 100   & 2.65 & 6.53 & 3.77e-03 & 3.79e-06 \\
             &          & 1000  & 2.65 & 3.43 & 1.54e-03 & 3.79e-06 \\
\addlinespace
$128\times128$ & 297478 & 1   & 2.85 & 10.88 & 1.55e-03 & 4.77e-07 \\
             &          & 10    & 2.90 & 8.57 & 1.49e-03 & 4.77e-07 \\
             &          & 100   & 2.70 & 7.11 & 1.24e-03 & 4.77e-07 \\
             &          & 1000  & 2.70 & 3.89 & 7.05e-04 & 4.88e-07 \\
\bottomrule
\end{tabular}
\end{table}

\subsubsection{Time-step robustness}

We next investigate the dependence of the solver on the time step $\Delta t$.  
The time step is varied from $5.0\times10^{-3}$ to $5.0\times10^{-4}$, while 
the final time is fixed at $T = 0.05$.  The mesh consists of $16\times16$ cells, 
the augmentation parameter is $\gamma = 100$, and the viscosity is $\nu = 0.1$.  
Table~\ref{tab:dt} gives the total number of time steps, the average Newton iterations per time step, the average FGMRES iterations per Newton step, and the final-time $L^2$ norms of the divergence and unit-length residuals.

\begin{table}[htbp]
\centering
\caption{Robustness with respect to the time step.  
         Mesh: $16\times16$, $\gamma = 100$, $\nu = 0.1$, $T = 0.05$.}
\label{tab:dt}
\begin{tabular}{cccccc}
\toprule
$\Delta t$ & Steps & Newton it & KSP it & $\|\nabla\cdot\mathbf{u}\|_{L^2}$ & $\||\mathbf{n}|^2-1\|_{L^2}$ \\
\midrule
0.005  & 10  & 3.40 & 5.62 & 7.68e-03 & 1.37e-04 \\
0.0025 & 20  & 2.45 & 5.41 & 7.48e-03 & 1.37e-04 \\
0.001  & 50  & 2.20 & 5.65 & 7.36e-03 & 1.36e-04 \\
0.0005 & 100 & 2.23 & 5.61 & 7.32e-03 & 1.36e-04 \\
\bottomrule
\end{tabular}
\end{table}
The average number of FGMRES iterations per Newton step ranges from $5.4$ to $5.7$ across the time steps considered, showing no systematic variation with $\Delta t$. The divergence residual
$\|\nabla\cdot\mathbf{u}\|_{L^2}$ and the director-length residual
$\||\mathbf{n}|^2-1\|_{L^2}$ remain essentially unchanged as
$\Delta t$ varies, indicating that these quantities are controlled
primarily by the spatial discretization rather than by the time step. The nearly constant outer count indicates time-step robustness over the single decade examined.

\subsubsection{Viscosity dependence}

The influence of the fluid viscosity $\nu$ is considered next. For fixed characteristic velocity and length scales, changing $\nu$ changes the Reynolds number. The viscosity is varied over three orders of magnitude: $\nu \in \{1.0,\,0.1,\,0.01,\,0.001\}$.  The mesh is fixed at $16\times16$ cells, $\gamma = 100$, $\Delta t = 10^{-3}$, and $T = 0.02$.  Table~\ref{tab:viscosity} summarises the results.

\begin{table}[htbp]
\centering
\caption{Robustness with respect to viscosity.  Mesh: $16\times16$, $\gamma = 100$, $\Delta t = 10^{-3}$, $T = 0.02$.}
\label{tab:viscosity}
\begin{tabular}{ccccccc}
\toprule
$\nu$ & Newton it & KSP it & $\|\nabla\cdot\mathbf{u}\|_{L^2}$ & $\||\mathbf{n}|^2-1\|_{L^2}$ & CPU (s) \\
\midrule
1.0   & 2.55 & 5.94 & 8.55e-03 & 2.31e-04 & 2.18 \\
0.1   & 2.50 & 5.82 & 1.90e-02 & 2.32e-04 & 1.99  \\
0.01  & 2.45 & 5.80 & 1.41e-02 & 2.33e-04 & 1.82 \\
0.001 & 2.45 & 5.80 & 6.03e-03 & 2.34e-04 & 1.80 \\
\bottomrule
\end{tabular}
\end{table}
The Newton iteration count is nearly constant across the whole range. The preconditioned FGMRES iteration count remains around six and shows no systematic growth as $\nu$ decreases, 
indicating that the outer Krylov convergence is largely independent of viscosity. The divergence residual varies mildly, while the director-length residual stays almost unchanged. 
These results support the conclusion that the Schur-complement approximation captures the dominant constraint contribution over the tested three-order-of-magnitude viscosity range.

\subsection{Flow and director evolution in a disk}

\begin{figure}[!ht]
  \centering
  \begin{subfigure}{0.48\linewidth}
    \centering
    \includegraphics[width=\linewidth]{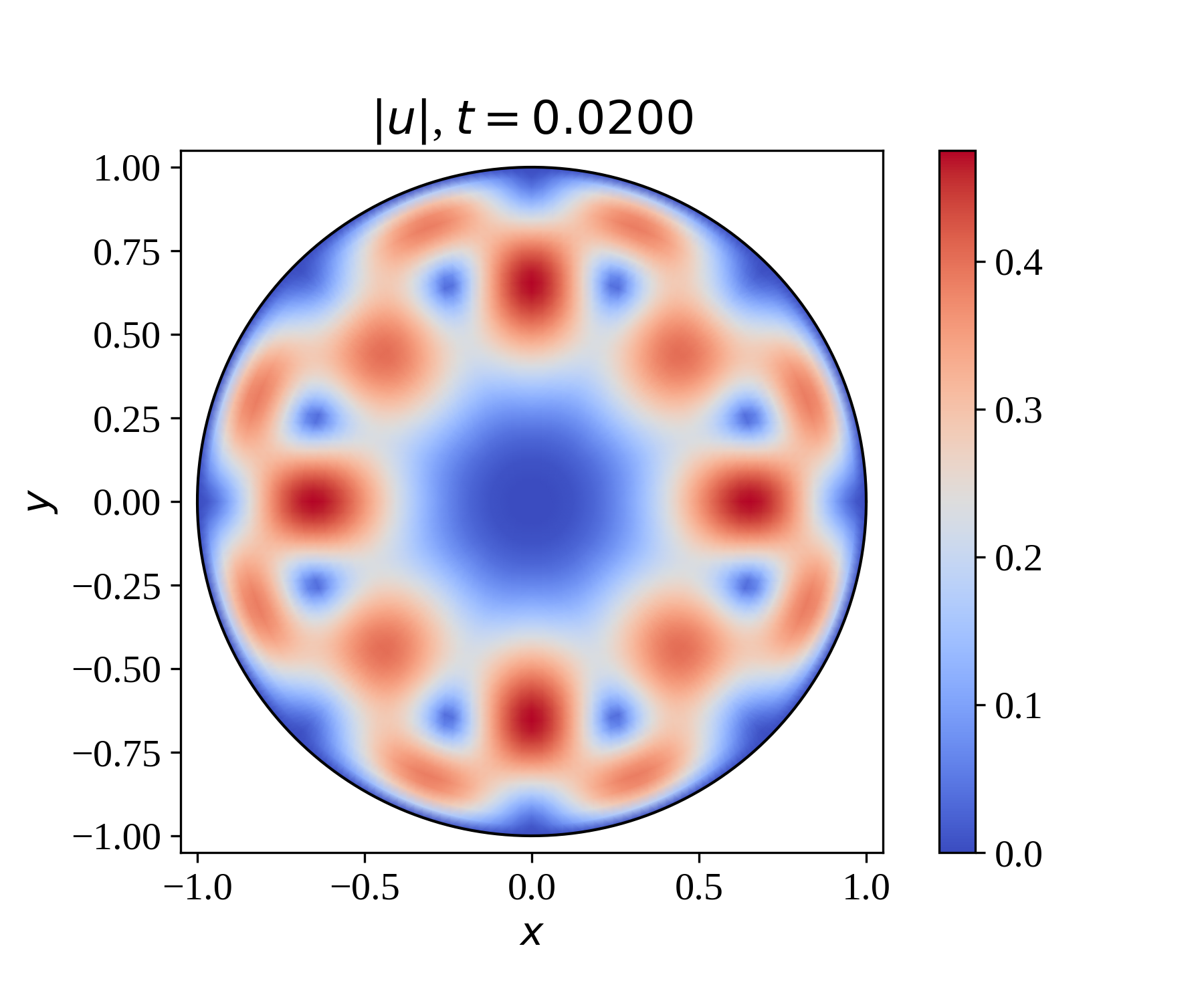}
    \caption{$t=0.02$}
  \end{subfigure}
  \hfill
  \begin{subfigure}{0.48\linewidth}
    \centering
    \includegraphics[width=\linewidth]{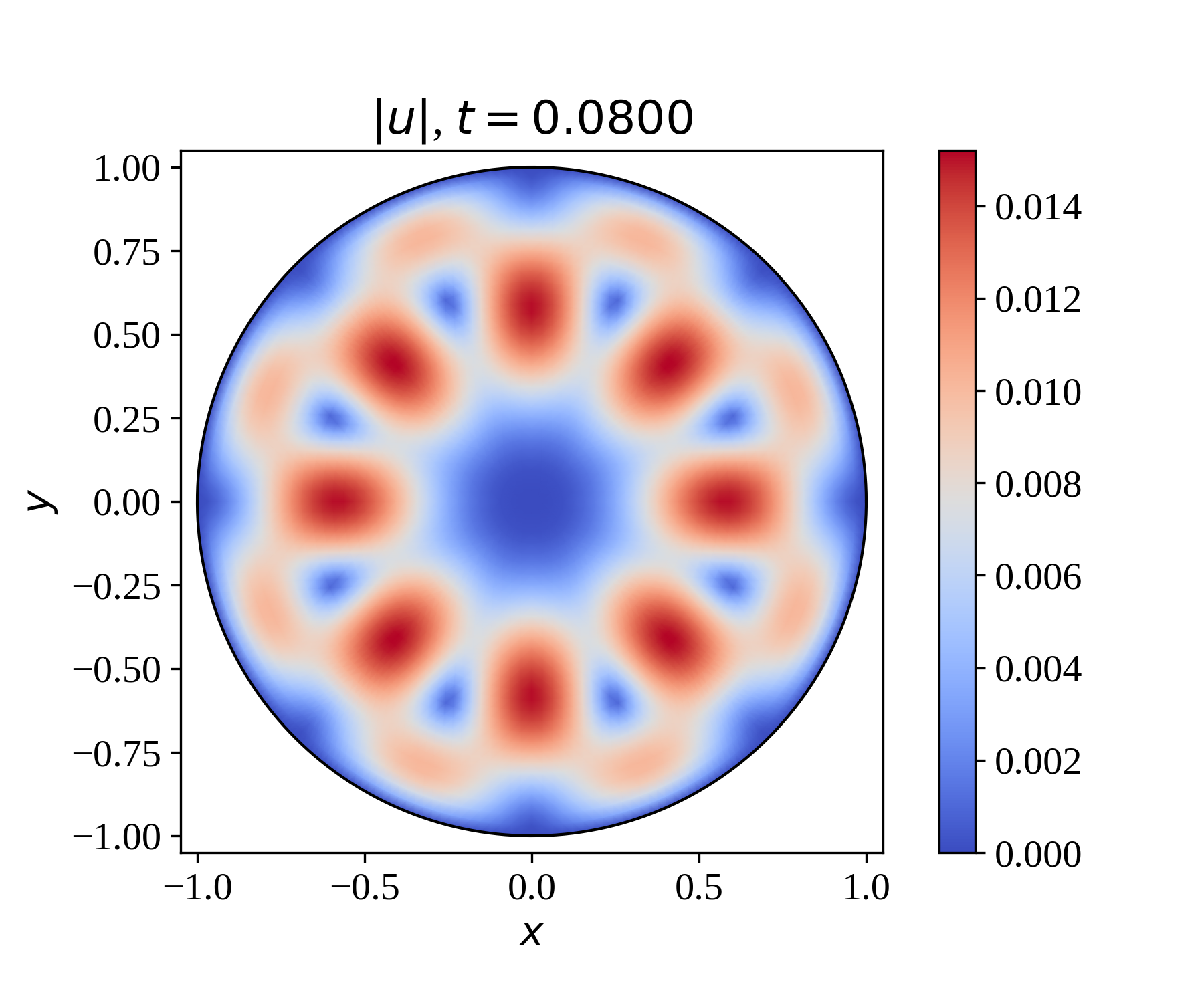}
    \caption{$t=0.08$}
  \end{subfigure}

  \vspace{0.5em}

  \begin{subfigure}{0.48\linewidth}
    \centering
    \includegraphics[width=\linewidth]{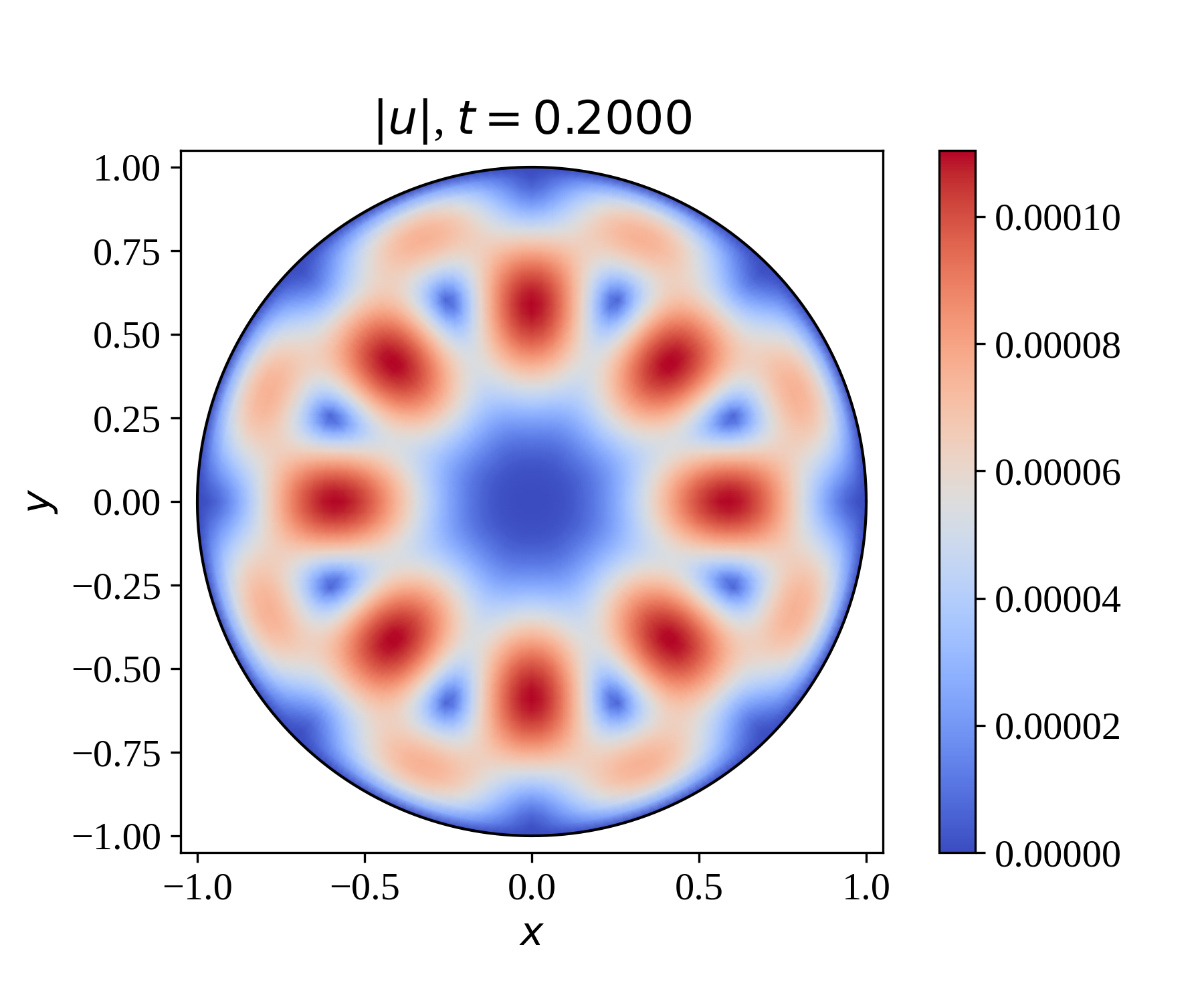}
    \caption{$t=0.2$}
  \end{subfigure}
  \hfill
  \begin{subfigure}{0.48\linewidth}
    \centering
    \includegraphics[width=\linewidth]{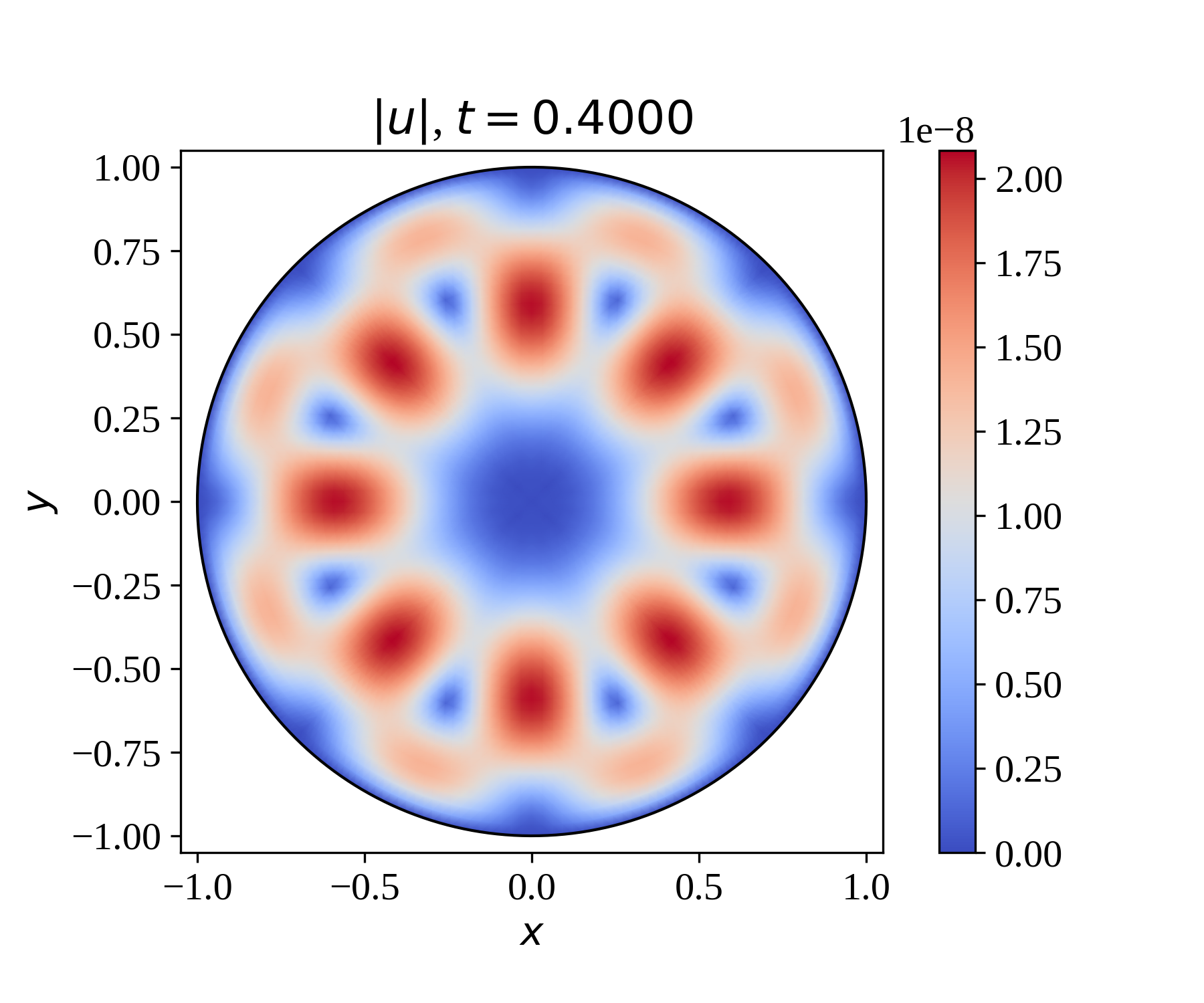}
    \caption{$t=0.4$}
  \end{subfigure}

  \caption{Velocity magnitude $|\mathbf{u}|$ at $t = 0.02,\,0.08,\,0.2,\,0.4$ for the
    unit disk problem.}
  \label{fig:disk_speed}
\end{figure}
To illustrate applicability on a curved geometry and a complex initial condition, we consider the test of Cao and Yi~\cite{CaoYi2025}. The computational domain is the unit disk $\Omega=\{x^2+y^2<1\}$ with no-slip boundary conditions for the velocity and homogeneous Neumann conditions for the director. The initial conditions are

\begin{equation}
\mathbf{u}_0 = \mathbf{0},\qquad 
\mathbf{n}_0 = (\sin a,\;\cos a)^{\!\top},\quad 
a = 4\pi\bigl(x^4 - y^4\bigr)^2,
\end{equation}
and the physical parameters are $\nu=2.0$, $K=1.0$, $\mu=1.0$. The augmentation parameters are $\gamma_u=\gamma_n=10.0$, the time step is $\Delta t = 1.0\times10^{-4}$, and the simulation is run until $T=0.4$. Spatial discretization uses $\mathbb{P}_2$ elements for $\mathbf{u}$ and $\mathbf{n}$, and $\mathbb{P}_1$ elements for $p$ and $q$ on a mesh obtained by five uniform refinements of a triangulation of the unit disk, resulting in approximately $2.1\times10^5$ degrees of freedom. The nonlinear system is solved with Newton's method, and the linearized saddle-point system is solved using FGMRES preconditioned by the augmented Lagrangian Schur complement approximations derived in Section~\ref{sec:Pre}.

\begin{figure}[!ht]
  \centering 
  \begin{subfigure}{0.48\linewidth}
    \centering
    \includegraphics[width=\linewidth]{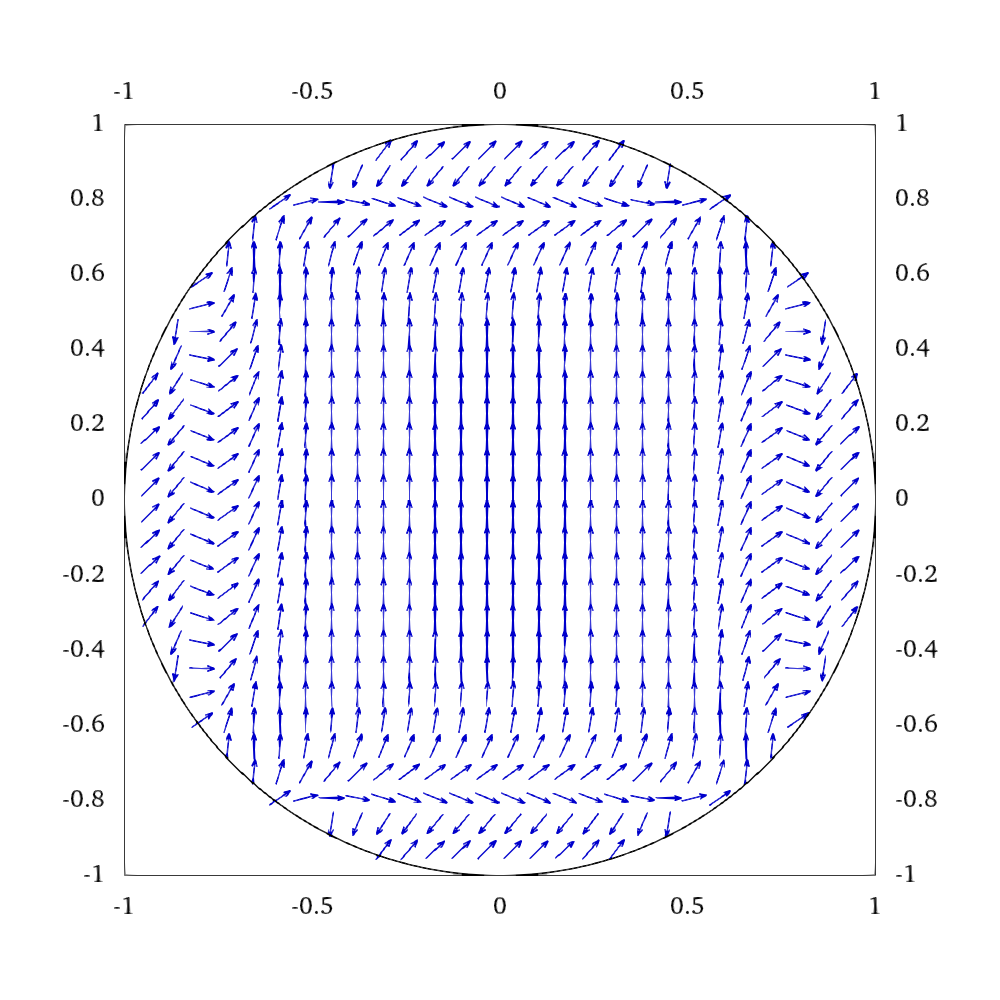}
    \caption{$t=0$}
  \end{subfigure}
  \hfill
  \begin{subfigure}{0.48\linewidth}
    \centering
    \includegraphics[width=\linewidth]{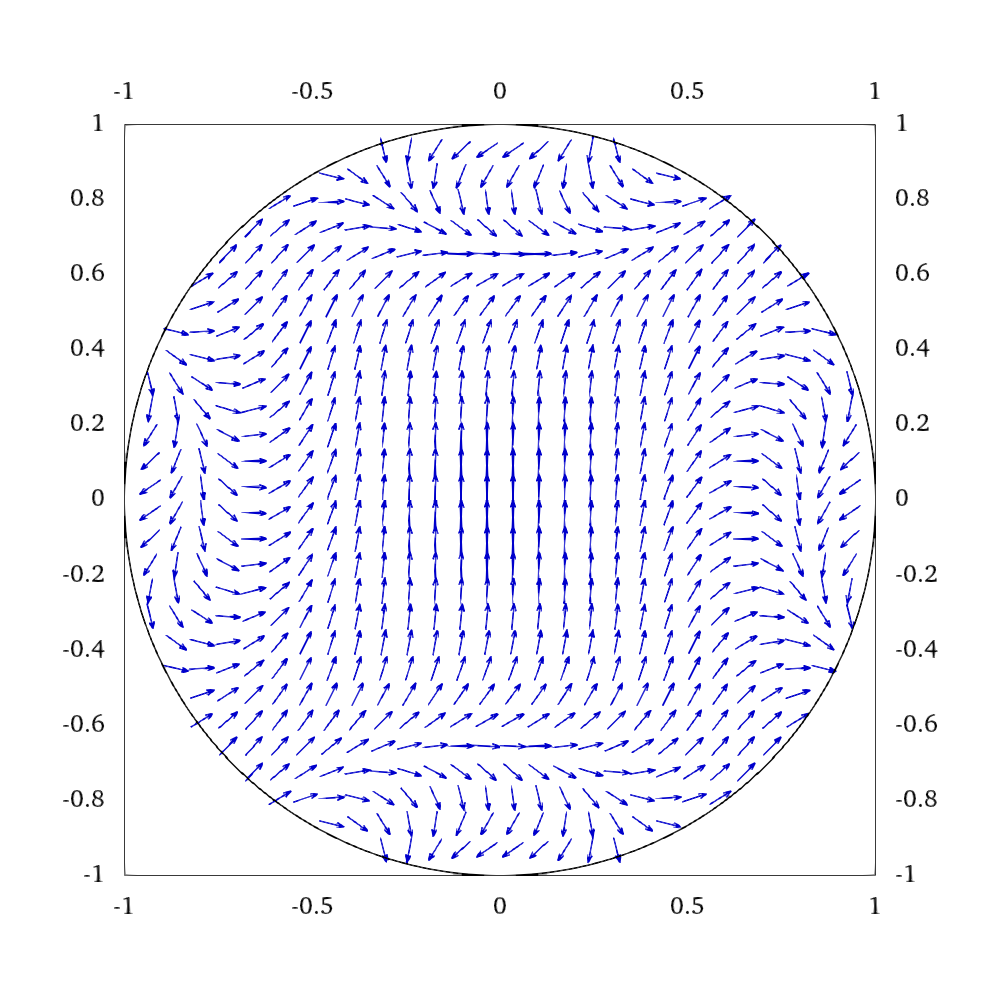}
    \caption{$t=0.02$}
  \end{subfigure}
  \vspace{0.5em}
  \begin{subfigure}{0.48\linewidth}
    \centering
    \includegraphics[width=\linewidth]{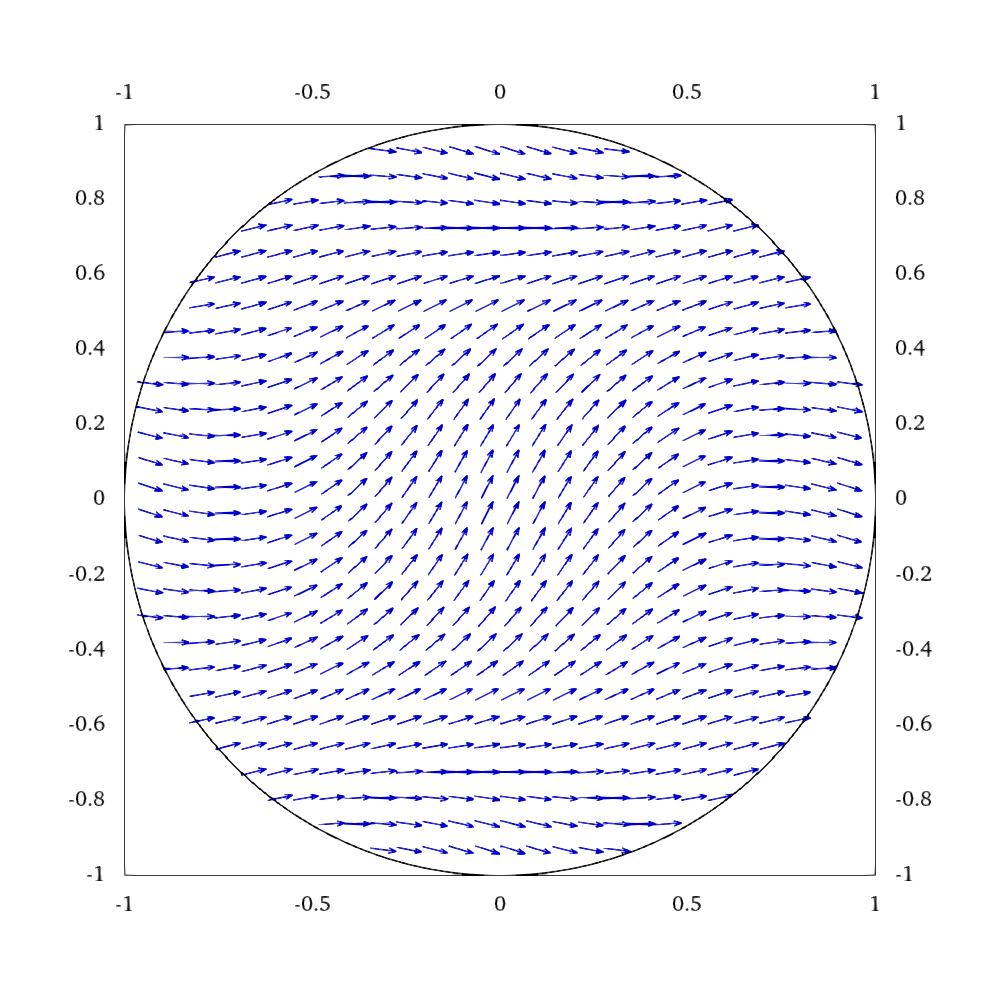}
    \caption{$t=0.08$}
  \end{subfigure}
  \hfill
  \begin{subfigure}{0.48\linewidth}
    \centering
    \includegraphics[width=\linewidth]{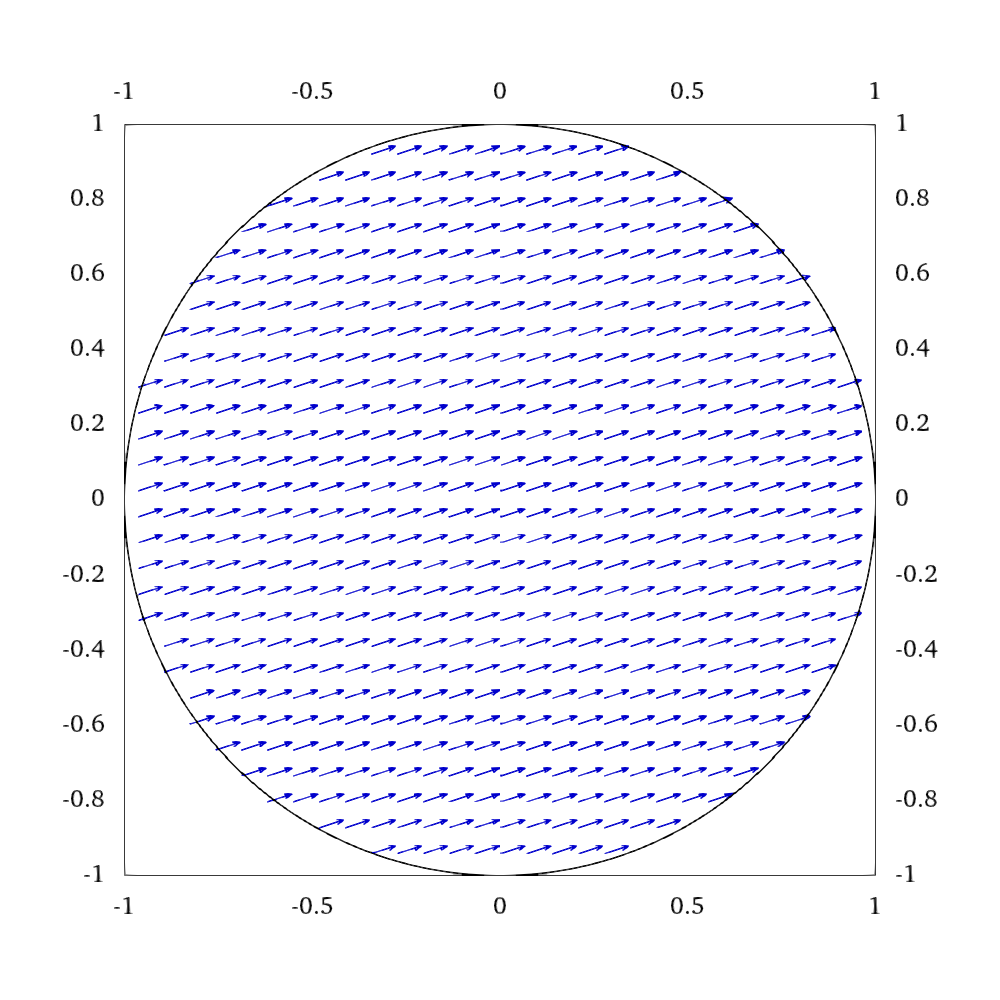}
    \caption{$t=0.4$}
  \end{subfigure}
  \caption{Director field $\mathbf{n}$ at $t=0$, $t=0.02$, $t=0.08$ and $t=0.4$ for the unit disk problem.}
  \label{fig:disk_director}
\end{figure}

Figures~\ref{fig:disk_speed} and~\ref{fig:disk_director} show that the velocity field rapidly relaxes while the director field evolves toward a smoother orientation pattern on the curved domain.  This test is qualitative, but it is useful for checking that the same block preconditioner can be applied without modification on a non-Cartesian mesh and for configurations with substantial director gradients.

\subsection{Two-defect configuration under a rotational flow}

As a final example, we consider a configuration containing two
near-singular regions in the director field, proposed by Badia
et al.~\cite{Badia2011}. This example involves an initial director field with large gradients and a strong rotational velocity that induces significant director reorientation, resulting in nontrivial coupling between flow and orientation throughout the simulation.  It serves to demonstrate the qualitative behavior of the augmented Lagrangian scheme in a regime close to physically relevant simulations.

The computational domain is $\Omega = (-1,1)^2$.  The initial director
field is constructed from the auxiliary vector field
\begin{equation}
\widetilde{\mathbf{n}}_0(x,y) =
\begin{pmatrix}
x^2 + y^2 - a^2 \\[2pt]
2a y
\end{pmatrix},
\qquad a = 0.5,
\end{equation}
followed away from the defect cores by the formal normalization
This defines two defects initially located at $(\pm a,0)$, where $|\widetilde{\mathbf{n}}_0| = 0$. In the finite element discretisation, to avoid division by zero during nodal interpolation at the defect cores, we introduce a small numerical parameter $\delta_{\text{init}} = 10^{-4}$ and define the discrete initial condition as
\[
\mathbf{n}_{0h} = \mathcal{I}_h \left( \frac{\widetilde{\mathbf{n}}_0}{\sqrt{|\widetilde{\mathbf{n}}_0|^2 + \delta_{\text{init}}^2}} \right).
\]
We emphasise that $\delta_{\text{init}}$ is employed solely for the initialisation step; no regularisation is applied during the time evolution, which is governed by the saddle-point model \eqref{eq:EL-saddle}. 
The initial velocity is prescribed as the rigid rotation
\[
\mathbf u_0=\omega(-y,x)^\top,\qquad \omega=50.
\]
The boundary conditions are no-slip 
for $\mathbf{u}$ and homogeneous Neumann for $\mathbf{n}$.

The physical parameters are $\nu = 1.0$, $K = 1.0$, $\mu = 1.0$, and the augmentation parameters are $\gamma_u = \gamma_n = 100$.  The mesh is a uniform triangular grid with $64 \times 64$ cells, and the time step is $\Delta t = 10^{-3}$. Spatial discretization uses $\mathbb{P}_2$ elements for $\mathbf{u}$ and $\mathbf{n}$, and $\mathbb{P}_1$ elements for $p$ and $q$, as in the previous sections.  The nonlinear system is solved with Newton's method, and the linearized saddle-point systems are solved with FGMRES preconditioned by the augmented Lagrangian Schur complement approximations described in Section~\ref{sec:Pre}.

\begin{figure}[!ht]
  \centering
  \begin{subfigure}{0.48\linewidth}
    \centering
    \includegraphics[width=\linewidth]{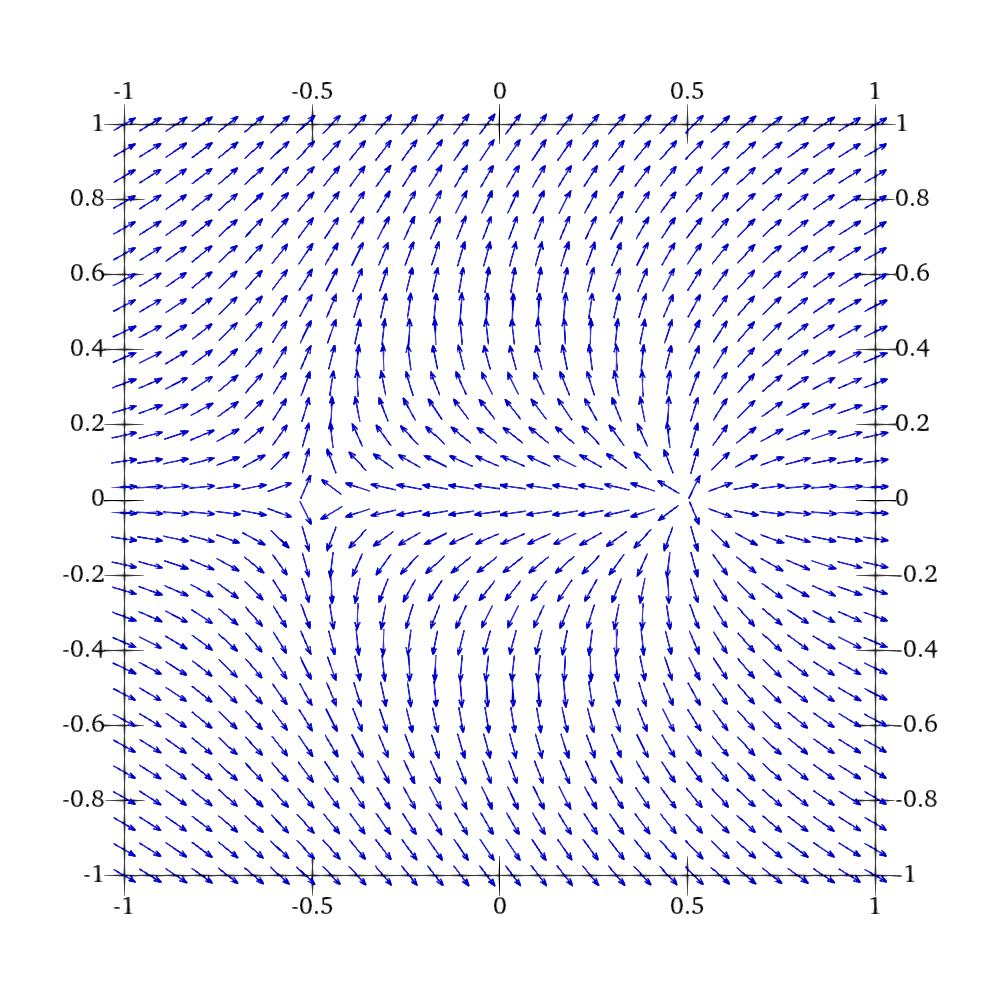}
    \caption{$t=0$}
  \end{subfigure}
  \hfill
  \begin{subfigure}{0.48\linewidth}
    \centering
    \includegraphics[width=\linewidth]{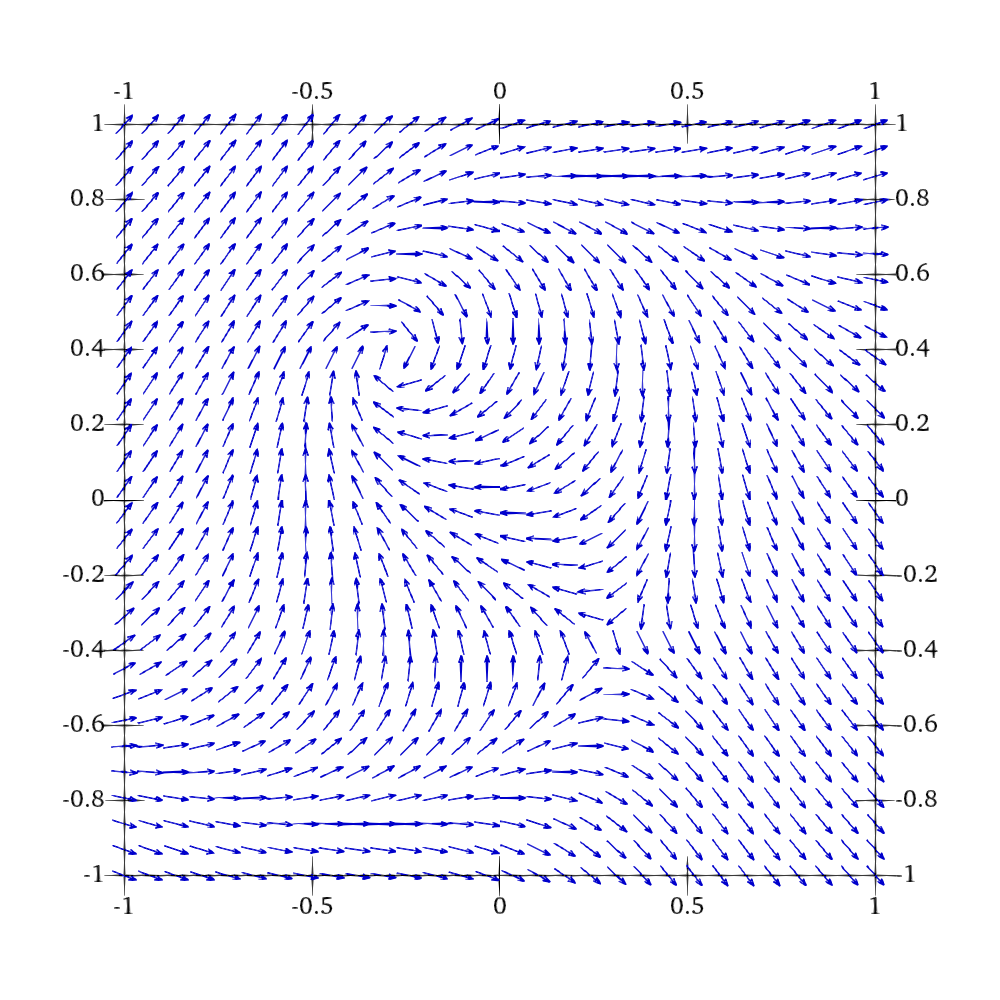}
    \caption{$t=0.05$}
  \end{subfigure}

  \vspace{0.5em}

  \begin{subfigure}{0.48\linewidth}
    \centering
    \includegraphics[width=\linewidth]{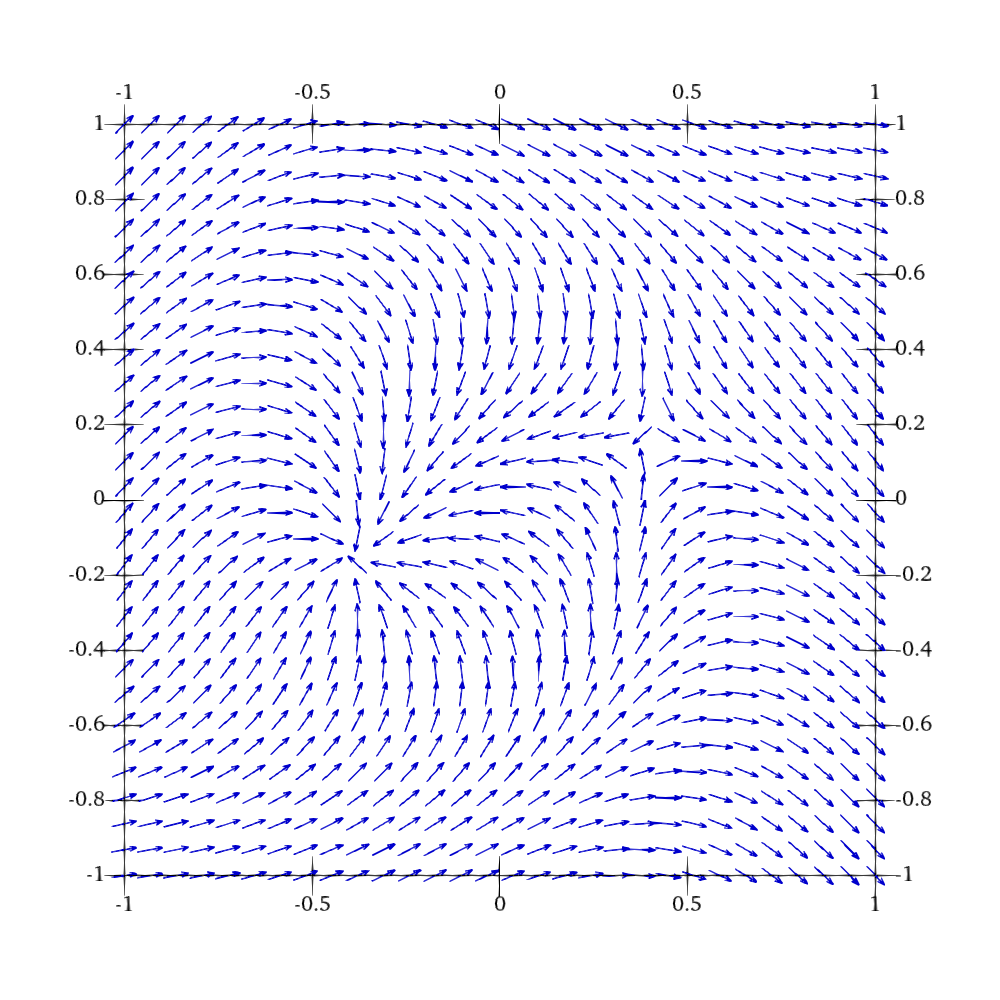}
    \caption{$t=0.10$}
  \end{subfigure}
  \hfill
  \begin{subfigure}{0.48\linewidth}
    \centering
    \includegraphics[width=\linewidth]{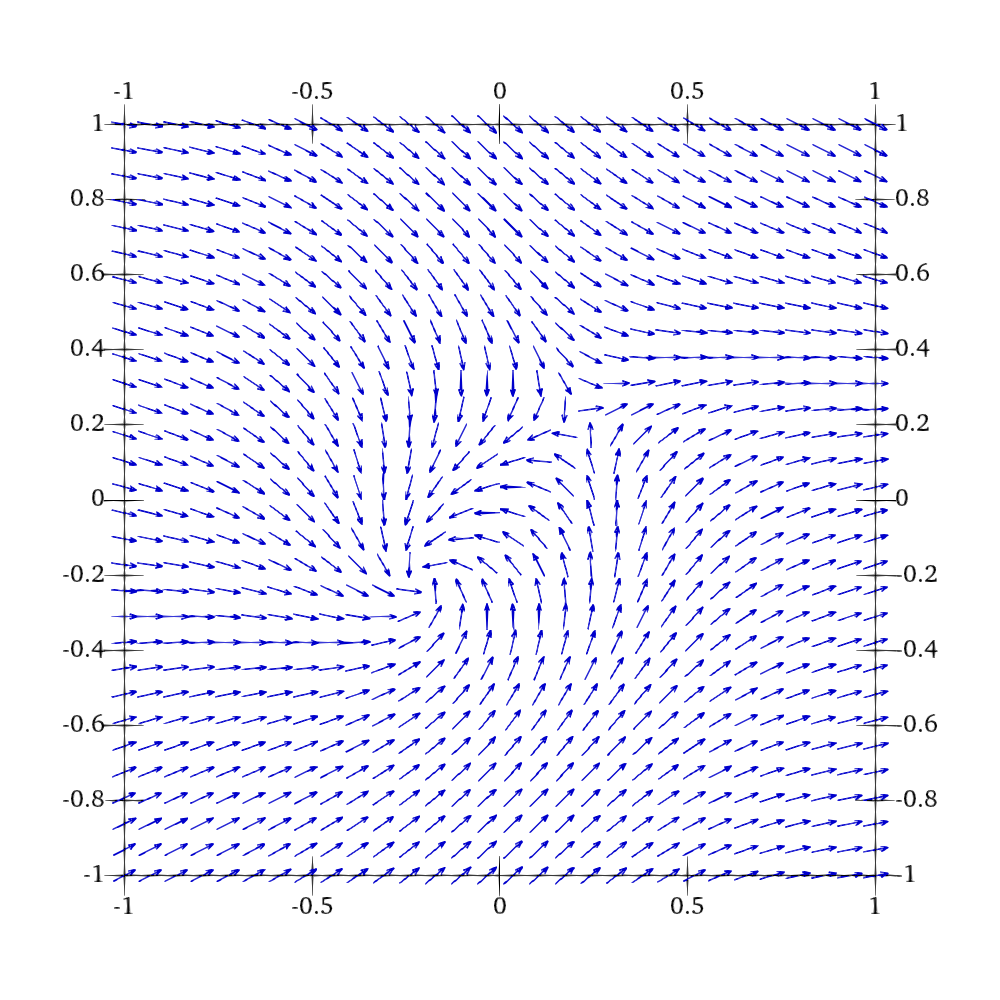}
    \caption{$t=0.25$}
  \end{subfigure}

  \caption{Director line field at $t = 0,\,0.05,\,0.10,\,0.25$ for the rotating
    two-defect problem. The two high-gradient regions near $(\pm0.5,0)$ are
    stretched and deformed by the rotational flow, and persist as stable
    non-uniform structures. Results obtained with the augmented Lagrangian
    preconditioner ($\gamma_u=\gamma_n=100$) on a $64\times64$ mesh.}
  \label{fig:badia_twodefect}
\end{figure}

Figure~\ref{fig:badia_twodefect} illustrates the spatial distributions of the director field at several representative time instants. Driven by the rotational flow, the two high-gradient regions are advected around the computational domain, and the surrounding orientation pattern undergoes substantial deformation. As the initial kinetic energy dissipates progressively, elastic relaxation gradually dominates the system evolution, and the director field eventually settles into a distorted quasi-steady configuration. The Newton--Krylov solver converges throughout the reported run. Apart from the explicit core-value convention used to define the initial finite element field, no smoothing is applied during the time evolution.

\section{Conclusions and perspectives}
\label{sec:conclusions}
We have developed an augmented Lagrangian block preconditioner for the Newton systems arising from a multiplier formulation of the simplified Ericksen--Leslie model. A block-diagonal approximation of the coupled velocity--director block leads to separate pressure and director-multiplier Schur approximations. The pressure scaling follows from the generalized-Stokes principal part and a Woodbury identity, whereas the multiplier scaling is motivated by the zero-order constraint operator and the reaction--diffusion principal part of the director equation.

The numerical results show the expected accuracy for the primary variables, nearly mesh-independent outer FGMRES counts over the tested meshes, stable behavior under time-step and viscosity variation, and improved iteration counts as the augmentation parameters increase. The energy benchmark is consistent with the continuous dissipation law. These conclusions are deliberately limited to the reported parameter ranges and to outer iteration counts.

A current limitation of this work is its restriction to the incompressible regime, which is the standard assumption for the simplified Ericksen--Leslie model. However, many nematic materials and polyatomic gases exhibit significant density variations. In a recent kinetic derivation for rarefied calamitic gases, Farrell, Russo, and Zerbinati \cite{farrell2024kinetic} derived an inviscid compressible variant of the Leslie--Ericksen equations, featuring a pressure-dependent Oseen--Frank energy functional. The non-trivial coupling between the fluid density and the nematic ordering introduced by this model presents a natural and challenging extension for the present preconditioning framework. Future work will focus on adapting the block preconditioning strategy to the compressible setting, which will require a reassessment of the saddle-point structure and the associated Schur complement approximations.


\section*{Acknowledgments}
This work is partially supported by the Innovation Research Foundation of National University of Defense Technology, the National Natural Science Foundation of China (No. 12371374), and the Youth Elite Scientists Sponsorship Program by CAST.

\section*{Code Availability Statement}
The source code used to generate the numerical examples presented in this manuscript is openly available in an repository at \url{https://github.com/LiYoca/al-preconditioner-ericksen-leslie}.

\bibliographystyle{elsarticle-num-names}
\bibliography{main}

@article{cao2026,
    author = {Ruonan, Cao and Nianyu Yi},
    title = {A linear, unconditionally stable, second order decoupled method for the {Ericksen-Leslie} model with {SAV} approach},
    journal = {Computers \& Mathematics with Applications},
    volume = {204},
    pages = {52--70},
    year = {2026}
}

@article{ericksen1962ARMA,
    author = {Ericksen, J. L.},
    title = {Hydrostatic theory of liquid crystals},
    journal = {Archive for Rational Mechanics and Analysis},
    pages = {371--378},
    number = {9},
    year = {1962} 
}

@article{ericksen1961,
    author = {Ericksen, J. L.},
    title = {Conservation Laws for Liquid Crystals},
    journal = {Journal of Rheology},
    pages = {23--34},
    number = {5},
    year = {1961}
}

@article{leslie1979,
    author = {Leslie, F. M.},
    title = {Theory of Flow Phenomena in Liquid Crystals},
    journal = {Advances in Liquid Crystals},
    pages = {1--81},
    volume = {4},
    year = {1979} 
}

@article{leslie1968ARMA,
    author = {Leslie, F. M.},
    title = {Some constitutive equations for liquid crystals},
    journal = {Archive for Rational Mechanics and Analysis},
    pages = {265--283},
    volume = {28},
    year = {1978} 
}

@article{lin1989,
    author = {Lin, F.},
    title = {Nonlinear theory of defects in nematic liquid crystals phase transition and flow phenomena},
    journal = {Commun. Pure Appl. Math.},
    pages = {789-–814},
    number = {42},
    year = {1989}
}

@article{saad1986,
  author = {Saad, Y. and Schultz, M. H.},
  title = {GMRES: a generalized minimal residual algorithm for solving nonsymmetric linear systems},
  journal = {SIAM J. Sci. Stat. Comput.},
  volume = {7},
  pages = {856--869},
  year = {1986}
}

@article{xia2021,
  author = {Xia, J. and Farrell, P. E. and Wechsung, F.},
  title = {Augmented Lagrangian preconditioners for the Oseen-{Frank} model of nematic and cholesteric liquid crystals},
  journal = {Numer. Math.},
  volume = {149},
  pages = {609--646},
  year = {2021}
}

@article{mardal2011,
  author = {Mardal, K.-A. and Winther, R.},
  title = {Preconditioning discretizations of systems of partial differential equations},
  journal = {Numer. Linear Algebra Appl.},
  volume = {18},
  pages = {1--40},
  year = {2011}
}

@article{cahouet1988,
  author = {Cahouet, J. and Chabard, J.-P.},
  title = {Some fast 3D finite element solvers for the generalized Stokes problem},
  journal = {Internat. J. Numer. Methods Fluids},
  volume = {8},
  pages = {869--895},
  year = {1988}
}

@book{elman2014,
  author = {Elman, H. C. and Silvester, D. J. and Wathen, A. J.},
  title = {Finite Elements and Fast Iterative Solvers: with Applications in Incompressible Fluid Dynamics},
  edition = {2nd},
  publisher = {Oxford University Press},
  address = {Oxford},
  year = {2014}
}

@article{schoberl1999,
  author = {Schöberl, J.},
  title = {Multigrid methods for a parameter dependent problem in primal variables},
  journal = {Numer. Math.},
  volume = {84},
  pages = {97--119},
  year = {1999}
}

@book{deGennesProst1993,
  author = {de Gennes, P. G. and Prost, J.},
  title = {The Physics of Liquid Crystals},
  edition = {2nd ed.},
  publisher = {Oxford University Press},
  year = {1993}
}

@book{Stewart2004,
  author = {Stewart, I. W.},
  title = {The Static and Dynamic Continuum Theory of Liquid Crystals},
  publisher = {Taylor \& Francis},
  year = {2004}
}

@article{LinLiu1995,
  author = {Lin, F. H. and Liu, C.},
  title = {Nonparabolic dissipative systems modeling the flow of liquid crystals},
  journal = {Comm. Pure Appl. Math.},
  volume = {48},
  pages = {501--537},
  year = {1995}
}

@article{Oseen1933,
  author = {Oseen, C. W.},
  title = {The theory of liquid crystals},
  journal = {Trans. Faraday Soc.},
  volume = {29},
  pages = {883--899},
  year = {1933}
}

@article{Frank1958,
  author = {Frank, F. C.},
  title = {On the theory of liquid crystals},
  journal = {Discuss. Faraday Soc.},
  volume = {25},
  pages = {19--28},
  year = {1958}
}

@article{Badia2011,
  author = {Badia, Santiago and Guillén-González, Francisco and Gutiérrez-Santacreu, Juan Vicente},
  title = {Finite element approximation of nematic liquid crystal flows using a saddle-point structure},
  journal = {Journal of Computational Physics},
  volume = {230},
  pages = {1686--1706},
  year = {2011}
}

@article{CaoYi2025,
  author = {Cao, Ruonan and Yi, Nianyu},
  title = {Length preserving numerical schemes for the nematic liquid crystal flows},
  journal = {ESAIM: Mathematical Modelling and Numerical Analysis},
  volume = {59},
  pages = {3021--3040},
  year = {2025}
}

@article{Farrell2019,
  author    = {Farrell, Patrick E. and Mitchell, Lawrence and Wechsung, Florian},
  title     = {An augmented Lagrangian preconditioner for the {3D} stationary incompressible {Navier–Stokes} equations at high Reynolds number},
  journal   = {SIAM Journal on Scientific Computing},
  volume    = {41},
  number    = {5},
  pages     = {A3073--A3096},
  year      = {2019}
}

@article{Benzi2006,
  author    = {Benzi, Michele and Olshanskii, Maxim A.},
  title     = {An augmented Lagrangian-based approach to the {Oseen} problem},
  journal   = {SIAM Journal on Scientific Computing},
  volume    = {28},
  number    = {6},
  pages     = {2095--2113},
  year      = {2006}
}

@article{Alouges1997,
  author  = {Alouges, Fran{\c{c}}ois},
  title   = {A new algorithm for computing liquid crystal stable configurations: 
             the harmonic mapping case},
  journal = {SIAM Journal on Numerical Analysis},
  volume  = {34},
  number  = {5},
  pages   = {1708--1726},
  year    = {1997}
}

@article{farrell2024kinetic,
    author  = {P. E. Farrell and G. Russo and U. Zerbinati},
    year    = {2024},
    title   = {{Kinetic Derivation of an Inviscid Compressible Leslie--Ericksen Equation for Rarified Calamitic Gases}},
    journal = {Multiscale Modeling \& Simulation},
    volume  = {22},
    number  = {4},
    pages   = {1585--1607},
}

@article{Guti2017,
author = {Guti\'{e}rrez-Santacreu, Juan Vicente and Restelli, Marco},
title = {Inf-Sup Stable Finite Element Methods for the Landau--Lifshitz--Gilbert and Harmonic Map Heat Flow Equations},
journal = {SIAM Journal on Numerical Analysis},
volume = {55},
number = {6},
pages = {2565-2591},
year = {2017},
}

@article{Hu2009,
author = {Hu, Qiya and Tai, Xue-Cheng and Winther, Ragnar},
title = {A Saddle Point Approach to the Computation of Harmonic Maps},
journal = {SIAM Journal on Numerical Analysis},
volume = {47},
number = {2},
pages = {1500-1523},
year = {2009},
}

@article{Bacuta2006,
author = {Bacuta, Constantin},
title = {A Unified Approach for Uzawa Algorithms},
journal = {SIAM Journal on Numerical Analysis},
volume = {44},
number = {6},
pages = {2633-2649},
year = {2006},
}

@article{firedrake2016,
author = {Rathgeber, Florian and Ham, David A. and Mitchell, Lawrence and Lange, Michael and Luporini, Fabio and Mcrae, Andrew T. T. and Bercea, Gheorghe-Teodor and Markall, Graham R. and Kelly, Paul H. J.},
title = {Firedrake: Automating the Finite Element Method by Composing Abstractions},
year = {2016},
issue_date = {September 2017},
publisher = {Association for Computing Machinery},
address = {New York, NY, USA},
volume = {43},
number = {3},
issn = {0098-3500},
journal = {ACM Trans. Math. Softw.},
month = dec,
articleno = {24},
numpages = {27}
}

@article{farrell2021irksome,
author = {Farrell, Patrick E. and Kirby, Robert C. and Marchena-Men\'{e}ndez, Jorge},
title = {Irksome: Automating Runge–Kutta Time-stepping for Finite Element Methods},
year = {2021},
issue_date = {December 2021},
publisher = {Association for Computing Machinery},
address = {New York, NY, USA},
volume = {47},
number = {4},
issn = {0098-3500},
journal = {ACM Trans. Math. Softw.},
month = sep,
articleno = {30},
numpages = {26}
}

@article{2012Heister,
title={Efficient augmented Lagrangian-type preconditioning for the Oseen problem using Grad-Div stabilization},
author={Timo Heister and Gerd Rapin},
journal={International Journal for Numerical Methods in Fluids},
issue={1},
pages={118-134},
year={2013},
}
\end{document}